\documentclass[10pt]{amsart}
\usepackage[letterpaper, margin=1in]{geometry}
\usepackage{fancyhdr}
\usepackage{amsmath,amsthm,amssymb,bbm,hyperref,xcolor,qtree,pdflscape}
\usepackage{enumerate}
\usepackage{hyperref}
\usepackage{verbatim}
\usepackage{array}
\usepackage{mathrsfs}
\usepackage{todonotes}
\usepackage{pgfplots}
\usepackage{multicol}
\usepackage{subcaption}
\pgfplotsset{compat = 1.13}
\usepackage{tabu}
\usepackage{graphicx}
\usepackage{enumitem}
\usepackage{siunitx}
\usepackage{multirow}
\usepackage{systeme}
\graphicspath{ {./images/}  }
\theoremstyle{remark}

\newcommand\ii\item
\newcommand\ol\overline

\newcommand\half{\frac12}

\newcommand{\norm}[1]{\left\| #1 \right\|}
\newcommand\RR{\mathbb R}
\newcommand\ZZ{\mathbb Z}

\renewcommand\div{\mathrm{div} \,}
\newcommand\curl{\mathrm{curl} \,}

\newcommand\tr{\mathrm{tr} \,}

\newcommand{\myang}[1]{\left \langle #1 \right \rangle}
\newcommand{\myvert}[1]{\left \lvert #1 \right \rvert}
\newcommand{\myVert}[1]{\left \lVert #1 \right \rVert}
\newcommand{\paren}[1]{\left ( #1 \right )}

\newcommand\dd{\partial}

\newcommand{\R}{\mathbb{R}}
\newcommand{\bn}{\mathbf{n}}

\newtheorem{theorem}{Theorem}[section]

\newtheorem{lemma}[theorem]{Lemma}

\newtheorem{Remark}[theorem]{Remark}
\newtheorem{scheme}[theorem]{Scheme}
\numberwithin{equation}{section}
\newenvironment{myeq}{\begin{equation}\begin{aligned}}{\end{aligned}\end{equation}}
\newenvironment{myeq*}{\begin{equation*}\begin{aligned}}{\end{aligned}\end{equation*}}

\title[Quartic Q-tensors]{Finite element analysis of a nematic liquid crystal Landau-de Gennes model with quartic elastic terms}
\date{\today}

\author[J. Elafandi]{Jacob Elafandi}
\address[Jacob Elafandi]{\newline Department of Mathematics \newline UC Berkeley \newline Berkeley, CA, 94720, USA.}
\email[]{elafandi@berkeley.edu}
\author[F. Weber]{Franziska Weber}
\address[Franziska Weber]{\newline Department of Mathematics \newline UC Berkeley \newline Berkeley, CA, 94720, USA.}
\email[]{fweber@berkeley.edu}

\thanks{F.W. and J.E. were supported in part by NSF grants DMS-2438083 and DMS-2409989.}

\begin{document}

\begin{abstract}
Golovaty et al.~\cite{GOLOVATY_2020} present a $Q$-tensor model for liquid crystal dynamics which reduces to the well-known Oseen-Frank director field model in uniaxial states. This model has been shown to capture phase transitions in nematic liquid crystals. We study a closely related model and present an energy stable scheme for the corresponding gradient flow. We prove well-posedness of the scheme via the Leray-Schauder fixed point theorem and rigorously show the $\Gamma$-convergence of discrete minimizers as the mesh size approaches zero. In the numerical experiments, we successfully simulate isotropic-to-nematic phase transitions as expected. We also compare simulations for elastic energies with identical and with highly disparate elastic constants and observe different defect dynamics. A comparison with simulations for the standard Landau-de Gennes model with quadratic energy shows that the quartic model is able to capture richer dynamics than the quadratic model.
\end{abstract}
\maketitle

\section{Introduction}
Liquid crystal is an intermediate state of matter between the liquid and crystalline solid states that exhibits a state of partial order~\cite{Stewart2019}. There are several types of liquid crystal phases; one is the nematic phase. In a material exhibiting this phase, the molecules can be represented as elongated rods, and intermolecular forces are responsible for the alignment of the molecules in a particular direction. Such materials are referred to (through a slight abuse of terminology) as nematic liquid crystals. Because the molecules in a liquid crystal are aligned but retain the ability to flow, liquid crystals are attractive for many applications in engineering and daily life, such as displays, smart glasses, and soaps~\cite{Stewart2019,Castellano2005}.

\subsection{Mathematical models for liquid crystal dynamics}
A variety of mathematical models for the dynamics of liquid crystals have been proposed. Two of the most commonly used ones are the Oseen-Frank director field model~\cite{Stewart2019,Frank1958,Oseen1933} and the Landau-de Gennes $Q$-tensor model~\cite{deGennes1995,Virga1995}. In the director field model, the main orientation of the liquid crystal molecules is modeled as a normal vector field $\bn$, which points in the general direction of the liquid crystal molecules. The direction of $\bn$ at any point in space is determined by minimizing the Oseen-Frank energy
\begin{equation}
\label{eq:OF}
\begin{split}
 \mathcal{F}_{OF}(\bn) &= \int_\Omega \bigg(
\frac{K_1}2(\div \bn)^2 +
\frac{K_2}2((\curl \bn) \cdot \bn)^2 +
\frac{K_3}2 |(\curl \bn) \times \bn|^2 \\ 
&\qquad +
\frac{K_2+K_4}2(\tr(\nabla \bn)^2 - (\div \bn)^2)
\bigg) \, dx,
\end{split}
\end{equation}
where $K_1, K_2, K_3, K_4 \geq 0$ are material and temperature dependent elastic constants and $\Omega\subset \R^3$ is the domain in which the liquid crystal is contained. Assuming all derivatives exist and commute, the fourth term in the integral can be rewritten as
\begin{equation*}
\tr(\nabla \bn)^2 - (\div \bn)^2 = \div\left((\bn\cdot \nabla)\bn-\bn \, \div \bn\right),
\end{equation*}
and so its integral can be expressed as a surface integral via Gauss' divergence theorem. Therefore this term is a null Lagrangian and does not contribute to the Euler-Lagrange equations which describe the equilibrium, although it enters the boundary conditions and affects the total value of the energy~\cite{Stewart2019}. This model can only account for uniaxial configurations in liquid crystals.

A theory that can also account for biaxial configurations is the Landau-de Gennes Q-tensor theory.
Here one assumes that the direction which the liquid crystal molecules favor has a probability density at every point in space and time. The first nontrivial moment of this probability density is the second moment $M$. If the material is in an isotropic phase at a point, the second moment reduces to $\frac{1}{3} I$. The $Q$-tensor is then defined as $Q = M -\frac{1}{3}I$, a traceless and symmetric tensor corresponding to the deviation from the isotropic phase. We thus define
\begin{myeq*}
	\mathcal S = \left\{Q \in \RR^{3 \times 3}: Q^\top = Q, \tr Q = 0\right\},
\end{myeq*}
the space of symmetric and traceless matrices. In the standard Landau-de Gennes model, we can find the equilibrium configuration $Q:\Omega\to \mathcal{S}$ given a domain $\Omega$ and a boundary condition $\left. Q\right|_{\partial\Omega} =G:\partial\Omega\to \mathcal{S}$ by minimizing the Landau-de Gennes energy
\begin{myeq}\label{eq:FLdG}
	\mathcal{F}_{LdG}(Q) = \int_\Omega \sum_{i,j,k=1}^3 \paren{\frac{L_1}2 Q_{ij,k}^2 + \frac{L_2}2 Q_{ij,j} Q_{ik,k} + \frac{L_3}2 Q_{ik,j} Q_{ij,k} + W(Q)} \, dx,
\end{myeq}
where $Q_{ij}$ is the $ij$-th element of $Q$ and $Q_{ij,k} = \partial Q_{ij} / \partial x_k$. Here $L_1, L_2, L_3$ are material dependent elastic constants and $W$ is a bulk energy density given by
\begin{myeq}\label{eq:W}
	W(Q) = a \, \tr(Q^2) - \frac{2b}3 \, \tr(Q^3) + \frac c2(\tr(Q^2))^2.
\end{myeq}
This bulk energy is derived from a Taylor expansion about $Q = 0$~\cite{Mottram2014}. This expression has often been used in the literature, but other expressions -- such as a mean field approximation, as in~\cite{Schimming2021} -- are possible as well. We assume $c>0$ to ensure the energy has a lower bound. The signs of $a$ and $b$ vary with temperature, and we note in particular that $a$ is negative for sufficiently low temperatures. If $a<\frac{b^2}{27 c}$, which is what we will assume in this work, then $Q$-tensors in the minimal set $\mathcal{N}$ of $W$ are in a uniaxial nematic state that can be described in terms of the set 
\begin{equation*}
	\mathcal N = \left\{Q= s_0 \paren{\mathbf{n} \otimes \mathbf{n} - \frac 13 I}\,:\,\mathbf{n} \in \mathbb S^{2} \right\}.
\end{equation*}
Here $s_0$ depends on the coefficients $a,b,c$, specifically,
\begin{equation*}
s_0=\frac{b+\sqrt{b^2-24ac}}{4 c};
\end{equation*}
cf. ~\cite[Section II.A]{Mottram2014}. Researchers have studied the connections between the Oseen-Frank theory and the Landau-de Gennes $Q$-tensor model for several decades ~\cite{Berreman1984,VanderZwan1988,Dickmann1995,Longa1989}. In particular, Majumdar and Zarnescu~\cite{Majumdar2010} have studied this issue in the one-constant setting, i.e., when $K_1=K_2=K_3=K$. Recently, Golovaty et al.~\cite{GOLOVATY_2020} studied connections between the two theories for the full Oseen-Frank model with possibly variable constants. This is an important problem as many materials indeed exhibit material constants $K_1$, $K_2$ and $K_3$ with highly disparate magnitude.
Golovaty et al. suggested a higher-order generalization of the elastic energy in~\eqref{eq:FLdG} to include quartic terms:
\begin{myeq}\label{golovaty_energy}
	\mathcal{E}(Q) &= \int_\Omega \bigg(
	\frac {L_1}2 \myvert{\paren{\frac{s_0}3 I + Q} \div Q}^2 +
	\frac {L_2}2 \myvert{\paren{\frac{2s_0}3 I - Q} \div Q}^2 \\ &\quad +
	\frac {L_3}2 \myvert{\paren{\frac{s_0}3 I + Q} \curl Q}^2 +
	\frac {L_4}2 \myvert{\paren{\frac{2s_0}3 I - Q} \curl Q}^2 + W(Q) \bigg) \, dx
\end{myeq}
Here $W$ is given by~\eqref{eq:W} and $L_1,L_2,L_3,L_4 > 0$. Note that the terms in this model still have quadratic dependence on both $\nabla Q$ and $Q$.
This energy is derived by selecting quartic choices of elastic energy terms from a list of thirteen found in ~\cite{Longa1987}.
If $Q=s_0(\bn \otimes \bn -I/3)$, then $\mathcal E$ becomes the Oseen-Frank energy $\mathcal F_{OF}$ \eqref{eq:OF} for suitable choices of $s_0$, $L_i$ and $K_i$~\cite[Prop. 2.6]{GOLOVATY_2020}. Golovaty et al. establish a rigorous connection between $\mathcal E$ and $\mathcal F_{OF}$ via a $\Gamma$-convergence argument in the limit of vanishing nematic correlation length. The nematic correlation length is a measure of how far the orientation of molecules in a nematic liquid crystal is correlated. It represents the characteristic distance over which the nematic order parameter remains correlated. Essentially, it describes the extent of long-range orientational order in the material. Golovaty et al. define a parameter $\epsilon$ which    represents the ratio of the nematic correlation length to a characteristic lengthscale of the domain and nondimensionalize the energy to arrive at
\begin{myeq}\label{golovaty_energynondim}
	\mathcal{E}_\epsilon(Q) &= \int_\Omega \bigg(
	\frac {L_1}2 \myvert{\paren{\frac{s_0}3 I + Q} \div Q}^2 +
	\frac {L_2}2 \myvert{\paren{\frac{2s_0}3 I - Q} \div Q}^2 \\ &\quad +
	\frac {L_3}2 \myvert{\paren{\frac{s_0}3 I + Q} \curl Q}^2 +
	\frac {L_4}2 \myvert{\paren{\frac{2s_0}3 I - Q} \curl Q}^2 + \frac{1}{\epsilon^2}W(Q) \bigg) \, dx.
	\end{myeq}
Then they use a $\Gamma$-convergence argument to send $\epsilon\to 0$. This can also be seen as a low temperature limit in some cases~\cite{Gartland2018}. 
The case of variable constants had been studied prior to the work of Golovaty et al.~\cite{GOLOVATY_2020}, but previous works were based mainly on the observation that one can obtain the Oseen-Frank energy density from the Landau-de Gennes energy density in the uniaxial case by adding an additional elastic term that is cubic in $Q$ and its derivatives, specifically,
\begin{equation*}
    \sum_{\ell,k,i,j=1}^3 Q_{\ell k}\partial_k Q_{ij}\partial_\ell Q_{ij}.
\end{equation*}
However, from the perspective of energy minimization, this term is problematic, since including it to the energy density results in the energy  being unbounded from below. The energy density sugested by Golovaty et al. does not have this issue. It is further motivated by~\cite{Longa1989} who suggested adding higher order terms to the Landau-de Gennes energy density in the presence of biaxial fluctuations.

In this paper, we consider the numerical analysis of a quartic model inspired by~\eqref{golovaty_energy}. Specifically, we consider the energy
\begin{equation}
\label{eq:Fkappa}
\begin{split}
\mathcal{F}(Q) &= \int_\Omega \bigg(
\frac {L_1}2 \myvert{\paren{\frac{s_0}3 I + Q} \div Q}^2 +
\frac {L_2}2 \myvert{\paren{\frac{2s_0}3 I - Q} \div Q}^2 + \frac {L_3}2 \myvert{\paren{\frac{s_0}3 I + Q} \curl Q}^2 \\ 
&\qquad  +
\frac {L_4}2 \myvert{\paren{\frac{2s_0}3 I - Q} \curl Q}^2+ \frac{L_5}{2}\left|Q\right|^2\left|\nabla Q\right|^2  + W(Q) \bigg) \, dx,
\end{split}
\end{equation}
where $L_5\geq 0$ is another material constant.  $L_5=0$ corresponds to Golovaty et al.'s original model~\eqref{golovaty_energy}. The term $|Q|^2|\nabla Q|^2$ is one of the possible choices of quartic terms in the list found in~\cite{Longa1987}. 
Moreover, it can be motivated from the Oseen-Frank energy~\eqref{eq:OF} in a similar way as Golovaty et al. derived the other four elastic terms  from the Oseen-Frank energy in the uniaxial case in~\cite[Prop. 2.6]{GOLOVATY_2020}. To see this, we first note that the Oseen-Frank energy can be written in the form
\begin{equation*}
\begin{split}
&\mathcal{F}_{OF}(\bn) = \int_{\Omega} \bigg(\frac{K_1}{2}|\nabla \bn|^2 + \frac{K_2-K_3}{2}(\curl\bn\cdot \bn)^2+\frac{K_3-K_1}{2}|\curl\bn|^2 \\
&\hphantom{\mathcal{F}_{OF}(\bn) = \int_{\Omega} \bigg(}
+\frac{K_2+K_4-K_1}{2}(\tr(\nabla\bn)^2-(\div\bn)^2) \bigg) dx;
\end{split}
\end{equation*}
see~\cite{Stewart2019}.
Then since $2|\nabla \bn|^2 = |\nabla(\bn\otimes\bn)|^2=|\nabla(\bn\otimes\bn-I/3)|^2$, we can rewrite the first term as
\begin{equation*}
\frac{K_1}{2} |\nabla \bn|^2 = \frac{3K_1}{8} \left|\bn\otimes\bn-\frac{I}{3}\right|^2 \left|\nabla\left(\bn\otimes\bn-\frac{I}{3}\right)\right|^2 = \frac{3K_1}{8s_0^4}\left|s_0\left(\bn\otimes\bn-\frac{I}{3}\right)\right|^2 \left|\nabla\left(s_0\left(\bn\otimes\bn-\frac{I}{3}\right)\right)\right|^2,
\end{equation*}
which is equal to $\frac{3K_1}{8s_0^4}|Q|^2|\nabla Q|^2$ in the case that $Q=s_0 (\bn\otimes\bn-I/3)$. The remaining terms with elastic constants $L_1, L_2, L_3, L_4$ can be justified in the same way as it was previously done in the work by Golovaty et al.~\cite[Section 2]{GOLOVATY_2020} (in fact the $|\nabla \bn|^2$ also appears there but is then rewritten to take a different form). This yields a justification for including the additional term. Thus, the $L_5$-term can be naturally motivated from the Oseen-Frank model in the uniaxial case. The model~\eqref{eq:Fkappa}   still has the same favorable properties as the one by Golovaty et al.\eqref{golovaty_energy}: Firstly, all the theoretical results established in~\cite{GOLOVATY_2020} carry over to the model with energy $\mathcal{F}$ with simple adaptions of the proofs. In particular, for a reasonable range of elastic constants, one can recover the Oseen-Frank model in the limit of vanishing correlation length. We can therefore take highly disparate values for the elastic constants in~\eqref{eq:Fkappa} and  simulate liquid crystals which favor one of the three deformations, splay, twist and bend. This allows, as we shall see in the simulations in Section~\ref{sec:res}, to capture the formation of non-spherically symmetric tactoids and nematic bubbles with defects. The additional term also respects all necessary symmetries. It is mentioned in~\cite{GOLOVATY_2020}, that their choice of quartic elastic energy terms may not be the only possible higher order extensions of the Landau-de Gennes model and we agree with this statement. Other terms in the list in~\cite{Longa1987} may lead to models with similar favorable properties. In fact, we believe that the work of this article may very well be extended to such similar models as long as the $L_5$-term is included. The $L_5$-term appears crucial to show $\Gamma$-convergence of the discrete energies $\mathcal{F}_h$ as $h \to 0$, where $h>0$ is a discretization parameter associated with our scheme, and in particular, to construct the recovery secquence, as we will see in Section~\ref{sec:gamma} of this article. 
This implies that if $\{Q_h^*\}_{h>0}$ is a sequence of global discrete minimizers of $\{\mathcal{F}_h\}_{h>0}$, then every cluster point of the sequence is a global minimizer of $\mathcal{F}$. In general, we do not know if the minima of $\mathcal{F}$ and its corresponding discrete energies are global and unique. If they are not, our $\Gamma$-convergence result implies that if $Q^*$ is an isolated local minimizer of $\mathcal{F}$, then there exists a sequence $Q_h^*$ converging to $Q^*$ with $Q_h^*$ a local minimizer of $\mathcal{F}_h$ for $h$ sufficiently small; cf. ~\cite[Theorem 5.1]{Braides2014}).

In order to find discrete minimizers, we consider an energy-stable finite element discretization of the gradient flow for the energy~\eqref{eq:Fkappa}.
To obtain the gradient flow, we compute the variational derivative of the energy and set it equal to the negative of the time derivative of $Q$:
\begin{equation}
\label{eq:gradflow}
\partial_t Q = -M\mathcal{P}\left(\frac{\delta \mathcal{F}(Q)}{\delta Q}\right),
\end{equation}
subject to boundary conditions $\left. Q\right|_{\partial\Omega} = G$.
Here $M>0$ is a constant, $\mathcal{P}$ is the projection onto symmetric traceless tensors given by 
\begin{myeq}
    \mathcal P(A) = \half(A + A^\top) - \frac{\tr A}3 I \quad \text{ for }\, A\in \R^{3\times 3},
\end{myeq}
and $\frac{\delta \mathcal{F}(Q)}{\delta Q}$ is given by
\begin{align*}
\dfrac{\delta \mathcal{F}(Q)}{\delta Q} &= L_1 \paren{- \nabla(S_1^2 \,\div Q) + S_1 \,\div Q (\div Q)^\top}+ L_2 \paren{\curl(S_1^2 \,\curl Q)^\top + S_1 \,\curl Q (\curl Q)^\top} \\
&- L_3 \paren{ \nabla(S_2^2 \,\div Q) + S_2\,\div Q(\div Q)^\top} + L_4 \paren{\curl(S_2^2 \,\curl Q)^\top - S_2 \,\curl Q (\curl Q)^\top} \\
&+ L_5\left(-\nabla(|Q|^2\nabla Q)+|\nabla Q|^2 Q\right)+ 2aQ - 2b (Q^2)+ 2c \,\tr (Q^2)Q,
\end{align*}
where we define $S_1 = \frac{s_0}{3}I+Q$, $S_2 = \frac{2s_0}{3}I-Q$, and $[\nabla(|Q|^2\nabla Q)]_{\ell m}=\sum_{i,j,k}\partial_i(Q_{jk}^2\partial_i Q_{\ell m})$.
Smooth solutions of the gradient flow~\eqref{eq:gradflow} satisfy the energy dissipation law
\begin{equation}\label{eq:energybalancecontinuous}
\frac{d}{dt}\mathcal{F}(Q) = -M\int_{\Omega} |\partial_t Q|^2 \, dx.
\end{equation}
Thus, starting from an arbitrary smooth initial condition, the energy decreases monotonically along the gradient flow 
\begin{equation*}
\mathcal{F}(Q(t))=\mathcal{F}(Q_0)-M\int_0^t \int_{\Omega}|\partial_s Q(s)|^2 \, dx \, ds,
\end{equation*}
where $Q_0=Q(0,\cdot)$ is the initial condition.
Heuristically, as time goes to infinity, the solution reaches an equilibrium point which is a local minimizer of $\mathcal{F}$. 
In this work, we design a fully discrete finite element scheme for~\eqref{eq:gradflow} that satisfies a discrete energy dissipation law. The discrete energy dissipation law mimics the energy law of the continuous setting~\eqref{eq:energybalancecontinuous}. This implies stability of the scheme independently of the mesh size $h$ and the time step size $\Delta t$ and a monotone decrease of the discrete energy with each time step. Therefore, the approximations converge to local minima of the corresponding discrete energies $\mathcal{F}_h$ as time goes to infinity. 
We show that the discrete energies are coercive and weakly lower semi-continuous on $H^1$ as the discretization parameters go to zero.
If $L_5>0$, we use this to show that the corresponding discrete energies $\mathcal{F}_h$ converge to the continuous energy in the sense of $\Gamma$-convergence as $h\to 0$.

While we are not able to show the same $\Gamma$-convergence result for the case $L_5=0$, which would correspond to the elastic energy density suggested by Golovaty et al.~\cite{GOLOVATY_2020}, our simulations with the proposed numerical method for the gradient flow show similar dynamics for $L_5$ small or zero. We present numerical results showing isotropic-to-nematic phase transitions as proposed in~\cite{golovaty2019phase}. The proposed numerical scheme for our model is able to capture all the tested phase transitions very well. In the numerical experiments, we observe the formation of defects and shrinking tactoids that have non-spherically symmetric boundaries, as observed in physical experiments, c.f.~\cite[Figure 1 and 2]{golovaty2019phase}. The dynamics vary depending on the ratio of the elastic constants. We also compare the simulations to those for the standard Landau-de Gennes quadratic one-constant approximation elastic energy and observe richer dynamics for the quartic model.
Additionally, with the exception of the $\Gamma$-convergence result, all stability and related results for the numerical scheme also hold in the case that $L_5=0$. 

\subsection{Related work}
To the best of our knowledge, apart from the simulations in~\cite{golovaty2019phase,Koizumi2023}, no computational works are available for the generalized Landau-de Gennes model with quartic energy density terms, and our work provides the first provably energy stable and convergent numerical scheme for~\eqref{eq:gradflow}. Related works on the classical Landau-de Gennes model with quadratic elastic energy include: ~\cite{Borthagaray2020,Borthagaray2021}, which discuss the numerical approximation of uniaxially constrained Q-tensors; ~\cite{gudibanda2022convergence,ZHAO2017803,Cai2017}, which design and analyze numerical methods for the gradient flow; and ~\cite{Schimming2021,Weber2023,Hirsch2023}, which construct numerical methods for Q-tensors subject to fluid flows and external forces, such as magnetic and electric fields. Our modification of the model~\eqref{golovaty_energy} allows for rigorous numerical analysis of  a Landau-de Gennes model with quartic elastic terms. We believe that the techniques that we develop to show $\Gamma$-convergence for this model are valuable in the sense that they may be used to show $\Gamma$-convergence for other related models and other energy stable numerical schemes. Due to the matrix-valued field $Q$ being potentially unbounded and not smooth, we introduce several approximation schemes in order to show the existence of a recovery sequence for the discrete energies. Furthermore, we observe in the numerical experiments that by allowing variable constants and higher order elastic energy terms, the dynamics of the gradient flow are richer and allows the formation of tactoid with boundaries displaying kinks, as it has been observed in experiments~\cite{golovaty2019phase}.

\subsection{Outline of this paper}
In Section~\ref{sec:not}, we introduce necessary notation that will be used in the rest of the paper. In Section~\ref{sec:num}, we introduce the numerical scheme and prove its solvability and stability. Section~\ref{sec:gamma} is dedicated to the proof of $\Gamma$-convergence of the discrete energies to the continuous one. Section \ref{sec:res} contains numerical convergence tests and simulations of isotropic-to-nematic phase transitions.

\section{Preliminaries}\label{sec:not}
We begin with a flurry of necessary definitions. For $A,B \in \RR^{3 \times 3}$, we use the Frobenius product and norm:
\begin{multicols}{2}
\begin{itemize}
    \item $A:B = \sum_{i,j=1}^3 A_{ij}B_{ij}$,
    \item $|A| = \sqrt{A : A}$.
\end{itemize}
\end{multicols}
We extend this definition for the gradients of tensor-valued functions:
\begin{multicols}{2}
\begin{itemize}
    \item $(\nabla A):(\nabla B) = \sum_{i,j,k=1}^3 A_{ij,k}B_{ij,k}$,
    \item $|\nabla A| = \sqrt{(\nabla A):(\nabla A)}$.
\end{itemize}
\end{multicols}
Let $A_1$, $A_2$, $A_3$ be the rows of $A:\Omega \to \mathbb R^{3 \times 3}$. As per \cite{GOLOVATY_2020}, we define the divergence of $A$ as
\begin{equation*}
    \div A = \sum_{j=1}^d (\div A_j) \textbf{e}_j
\end{equation*}
and its curl as the tensor map such that
\begin{equation*}
    (\curl A) v \equiv \curl (A^\top v) \quad \text{for all } v \in \mathbb R^3.
\end{equation*}
Thus the $j$-th entry of $\div A$ is $\div A_j$, and the $j$-th column of $\curl A$ is $\curl A_j$.

We also define the $L^2$ products and norms for scalar-valued functions $f,g:\Omega \to \RR$, vector-valued functions $u,v : \Omega \to \RR^d$, tensor-valued functions $A,B : \Omega \to \RR^{d \times d}$, and gradients of tensor-valued functions:
\begin{multicols}{2}
\begin{itemize}
    \item $\myang{f,g} = \int_\Omega fg \, dx$,
    \item $\myang{u,v} = \int_\Omega u \cdot v \, dx$,
    \item $\myang{A,B} = \int_\Omega A:B \, dx$,
    \item $\myang{\nabla A,\nabla B} = \int_\Omega (\nabla A):(\nabla B) \, dx$,
    \item $\myVert{f}^2 = \myang{f,f}$,
    \item $\myVert{u}^2 = \myang{u,u}$,
    \item $\myVert{A}^2 = \myang{A,A}$,
    \item $\myVert{\nabla A}^2 = \myang{\nabla A,\nabla A}$.
\end{itemize}
\end{multicols}
For $Q \in \mathcal S$, we define the symmetric tensors
\begin{myeq}
    S_1(Q) = \frac{s_0}3 I + Q, \qquad
    S_2(Q) = \frac{2s_0}3 I - Q,
\end{myeq}
which we use to define the functionals for the energy components:
\begin{myeq}
    \mathcal{F}_1(Q) &= \frac{L_1}2 \myVert{S_1 \, \div Q}^2, \quad
    \mathcal{F}_2(Q) = \frac{L_2}2 \myVert{S_1 \, \curl Q}^2, \quad
    \mathcal{F}_3(Q) = \frac{L_3}2 \myVert{S_2 \, \div Q}^2, \\
    \mathcal{F}_4(Q) &= \frac{L_4}2 \myVert{S_2 \, \curl Q}^2, \quad
    \mathcal{F}_5(Q) = \frac{L_5}2 \myVert{Q \, |\nabla Q|}^2, \quad
    \mathcal{F}_6(Q) = \int_\Omega W(Q) \, dx.
\end{myeq}
Then the energy \eqref{eq:Fkappa} can be rewritten as
\begin{myeq}
    \mathcal{F}(Q) &= \mathcal{F}_1(Q) + \mathcal{F}_2(Q) + \mathcal{F}_3(Q) + \mathcal{F}_4(Q) + \mathcal{F}_5(Q) + \mathcal{F}_6(Q).
\end{myeq}

In the vein of \cite{ZHAO2017803}, we denote the discrete value of a variable $(\cdot)(x,t)$ at time $t^n$ as $(\cdot)^n(x)$, and use the notation
\begin{equation*}
    (\cdot)^{n + \half} = \frac{(\cdot)^{n+1} + (\cdot)^n} 2.
\end{equation*}
For integer $k$, we also use the shorthand $\relax [k] = \{1, \ldots, k\}$.

\section{Discretization of the Q-tensor gradient flow}\label{sec:num}
Next, we introduce the numerical scheme for the gradient flow~\eqref{eq:gradflow} and prove its energy stability and solvability for each time step.
Given a family of quasi-uniform triangulations $\Omega_h$ of $\Omega$ with maximum element diameter $h$, we define a number of dependent quantities. To begin, let $L$ be the number of elements; let $N$ be the number of interior nodes; let $V$ be the total number of nodes, including boundary nodes; let $\{\Omega_\ell\}_{\ell=1}^L$ be the set of elements; let $\{P_i\}_{i=1}^{N}$ be the set of interior nodes; and let $\{P_i\}_{i=N+1}^V$ be the set of boundary nodes.
To deal with the boundary condition $G$, we introduce the following approximation of an extension of $G$ into the domain $\Omega$:
We define $G_h:\Omega \to \RR$ as the unique function which is linear on each $\Omega_\ell$ and for which
\begin{myeq}
    G_h(P_i) = \begin{cases}
        0 & i \in [N] \\
        G(P_i) & i \in [V] \setminus [N]
    \end{cases},
\end{myeq}
and define the space and associated set
\begin{myeq}
    \mathcal T_h^0 &= \{(\varphi:\Omega \to \mathcal S): \varphi \in C^0, \varphi \rvert_{\Omega_\ell} \text{ linear } \forall \ell, \varphi \rvert_{\dd \Omega} = 0\}, \\
    \mathcal T_h^G &= \{\varphi + G_h: \varphi \in \mathcal T_h^0\}.
\end{myeq}
Here $\mathcal{T}_h^0$ is the finite element space for homogeneous Dirichlet boundary conditions and $\mathcal{T}_h^G$ is a set to accommodate the boundary conditions for $Q$. $\mathcal{T}_h^G$ will be the set in which we look for an approximation of $Q$ and $\mathcal{T}_h^0$ will be the space for the test functions.
For each $i \in [V]$, define the ``hat'' function $\psi_i:\Omega \to \RR$ as the unique function which is continuous in $\Omega$ and linear on each $\Omega_\ell$ and equal to $\delta_{ij}$ at each $P_j$.
Let $E_{ij} = \textbf{e}_i \otimes \textbf{e}_j$, and define
\begin{myeq}\label{eq:Rs}
    R_1 = E_{11} - E_{33}, \quad R_2 = E_{22} - E_{33}, \quad R_3 = E_{12} + E_{21}, \quad R_4 = E_{13} + E_{31}, \quad R_5 = E_{23} + E_{32}.
\end{myeq}
Then
\begin{myeq}\label{eq:Bh}
    \mathcal B_h = \{R_\alpha \psi_i : \alpha \in [5], i \in [N]\}
\end{myeq}
is a basis for $\mathcal T_h^0$.

Labeling $S_1(Q_h)$ and $S_2(Q_h)$ as $S_{1,h}$ and $S_{2,h}$, respectively, we present a scheme which finds a finite element approximation of the solution of the gradient flow~\eqref{eq:gradflow}.

\begin{scheme}\label{fullscheme}
    If $Q_h^0 \in \mathcal T_h^G$ is given, then for all $n \geq 0$, we solve for $Q_h^{n+1} \in \mathcal T_h^G$ through:
    \begin{myeq}
    \myang{\frac{Q_h^{n+1} - Q_h^n}{\Delta t}, \varphi} &= -M H^{n+\half}(\varphi) \qquad \forall \varphi \in \mathcal B_h \text{ (and therefore } \forall \varphi \in \mathcal T_h^0)\\
    \text{where}&\\
    H^{n+\half}(\varphi) &= H_1^{n+\half}(\varphi) + H_2^{n+\half}(\varphi) + H_3^{n+\half}(\varphi) + H_4^{n+\half}(\varphi) + H_5^{n+\half}(\varphi) + H_6^{n+\half}(\varphi),\\
    \text{and}&\\
        H_1^{n + \half}(\varphi) &= L_1 \myang{ (S_{1,h} \,\div Q_h)^{n + \half}, S_{1,h}^{n + \half} \,\div \varphi + \varphi \,\div Q_h^{n + \half}}, \\
        H_2^{n + \half}(\varphi) &= L_2 \myang{ (S_{1,h} \,\curl Q_h)^{n + \half}, S_{1,h}^{n + \half} \,\curl \varphi + \varphi \,\curl Q_h^{n + \half}}, \\
        H_3^{n + \half}(\varphi) &= L_3 \myang{ (S_{2,h} \,\div Q_h)^{n + \half}, S_{2,h}^{n + \half} \,\div \varphi - \varphi \,\div Q_h^{n + \half}}, \\
        H_4^{n + \half}(\varphi) &= L_4 \myang{ (S_{2,h} \,\curl Q_h)^{n + \half}, S_{2,h}^{n + \half} \,\curl \varphi - \varphi \,\curl Q_h^{n + \half}}, \\
        H_5^{n + \half}(\varphi) &= L_5 \paren{\myang{(|\nabla Q_h|^2)^{n+\half} Q_h^{n+\half}, \varphi} + \myang{(|Q_h|^2)^{n+\half} \nabla Q_h^{n+\half}, \nabla \varphi} }, \\
        H_6^{n + \half}(\varphi) &= \myang{2aQ_h^{n + \half}
            - \frac{2b}3 (2(Q_h^2)^{n+\half} + Q_h^{n+1}Q_h^n)
            + 2c\paren{\tr (Q_h^2)}^{n+\half} Q_h^{n + \half}, \varphi}.
    \end{myeq}
\end{scheme}

This scheme is nonlinearly implicit, so solvability is not a priori clear. We will show in Theorem \ref{uniquesol} that a unique solution exists using a fixed point argument. Note that, if a solution $Q_h^{n+1}\in \mathcal{T}_h^G$ can be found, it is automatically traceless and symmetric by the definition of $\mathcal{T}_h^G$, c.f.~\eqref{eq:Rs} and~\eqref{eq:Bh}. First, however, we show that any approximation computed by Scheme~\ref{fullscheme} satisfies a discrete energy dissipation law corresponding to~\eqref{eq:energybalancecontinuous} which implies stability of the scheme.

\begin{theorem}\label{discrete_dissipation}
    Scheme \ref{fullscheme} satisfies the discrete energy dissipation law
    \begin{myeq}
        \frac{\mathcal{F}(Q_h^{n+1}) - \mathcal{F}(Q_h^n)}{\Delta t} = -\frac 1M \myVert{\frac{Q_h^{n+1} - Q_h^n}{\Delta t}}^2.
    \end{myeq}
\end{theorem}
\begin{proof}
    We begin by analyzing the $\mathcal{F}_1$ term:
    \begin{equation*}
    	\mathcal{F}_1(Q_h^{n+1}) - \mathcal{F}_1(Q_h^n) = \frac{L_1} 2 \paren{ \myVert{S_{1,h}^{n+1} \,\div Q_h^{n+1}}^2 - \myVert{S_{1,h}^n \,\div Q_h^n}^2 }.
    \end{equation*}
    Using the identity $|a|^2-|b|^2= (a+b)\cdot (a-b)$ for vectors $a$ and $b$, and  $a_1\cdot b_1 -a_2\cdot b_2 = \frac12 (a_1+a_2)\cdot(b_1-b_2)+\frac12 (a_1-a_2)\cdot (b_1+b_2)$ for matrices $a_1$, $a_2$ and vectors $b_1,b_2$,
    we can rewrite this as follows,
      \begin{myeq*}
    	&\mathcal{F}_1(Q_h^{n+1}) - \mathcal{F}_1(Q_h^n)
    		= \frac{L_1} 2 \myang{ S_{1,h}^{n+1} \,\div Q_h^{n+1} + S_{1,h}^{n} \,\div Q_h^{n}, S_{1,h}^{n+1} \,\div Q_h^{n+1} - S_{1,h}^{n} \,\div Q_h^{n} } \\
    		&= L_1 \myang{ \frac{S_{1,h}^{n+1} \,\div Q_h^{n+1} + S_{1,h}^{n} \,\div Q_h^{n}} 2, \frac{S_{1,h}^{n+1} + S_{1,h}^n} 2 \,\div(Q_h^{n+1} - Q_h^n) + (S_{1,h}^{n+1} - S_{1,h}^n) \,\div\paren{\frac{Q_h^{n+1} + Q_h^n} 2} }\\
    			&= L_1 \myang{ (S_{1,h} \,\div Q_h)^{n + \half} , S_{1,h}^{n + \half} \,\div(Q_h^{n+1}-Q_h^n) + (Q_h^{n+1} - Q_h^n) \,\div Q_h^{n + \half} } .
    	\end{myeq*}
    	  Because $Q^{n+1}_h-Q^n_h\in \mathcal{T}^0_h$ is a valid test function, we therefore have 
    	\begin{myeq*}
    		\mathcal{F}_1(Q_h^{n+1}) - \mathcal{F}_1(Q_h^n) 
    			&= L_1 \myang{ (S_{1,h} \,\div Q_h)^{n + \half} , S_{1,h}^{n + \half} \,\div(Q_h^{n+1}-Q_h^n) + (Q_h^{n+1} - Q_h^n) \,\div Q_h^{n + \half} } \\
    		&= H_1^{n + \half}(Q_h^{n+1} - Q_h^n).
    	\end{myeq*}
    The evaluation of the $\mathcal{F}_2$, $\mathcal{F}_3$, and $\mathcal{F}_4$ terms are analogous. For $\mathcal{F}_5$, we have, using again  $a_1\cdot b_1 -a_2\cdot b_2 = \frac12 (a_1+a_2)\cdot(b_1-b_2)+\frac12 (a_1-a_2)\cdot (b_1+b_2)$, and then $|a|^2-|b|^2 = (a+b):(a-b)$,
    \begin{myeq*}
        &\mathcal{F}_5(Q_h^{n+1}) - \mathcal{F}_5(Q_h^n) \\
        &= \frac{L_5}2 \int_\Omega \paren{ |Q_h^{n+1}|^2 |\nabla Q_h^{n+1}|^2 - |Q_h^n|^2 |\nabla Q_h^n|^2 } \, dx \\
        &= \frac{L_5}2 \int_\Omega \paren{\frac{|\nabla Q_h^{n+1}|^2+|\nabla Q_h^n|^2}2 (|Q_h^{n+1}|^2 - |Q_h^n|^2) + \frac{|Q_h^{n+1}|^2+|Q_h^n|^2}2 (|\nabla Q_h^{n+1}|^2 - |\nabla Q_h^n|^2)} \, dx \\
        &= L_5 \int_\Omega \paren{(|\nabla Q_h|^2)^{n+\half} Q_h^{n+\half} : (Q_h^{n+1} - Q_h^n) + (|Q_h|^2)^{n+\half} \nabla Q_h^{n+\half} : (\nabla Q_h^{n+1} - \nabla Q_h^n)} \, dx.
       \end{myeq*}
      Now, since $Q^{n+1}_h-Q^n_h\in \mathcal{T}^0_h$, we can write the last expression as
       \begin{myeq*}
        &\mathcal{F}_5(Q_h^{n+1}) - \mathcal{F}_5(Q_h^n) \\
         &= L_5 \int_\Omega \paren{(|\nabla Q_h|^2)^{n+\half} Q_h^{n+\half} : (Q_h^{n+1} - Q_h^n) + (|Q_h|^2)^{n+\half} \nabla Q_h^{n+\half} : (\nabla Q_h^{n+1} - \nabla Q_h^n)} \, dx	 \\
        &= H_5^{n+\half}(Q_h^{n+1} - Q_h^n).
    \end{myeq*}
    Lastly, we examine the $\mathcal{F}_6$ term. Using similar identities, we compute
    \begin{myeq*}
        \tr((Q_h^{n+1})^2) - \tr((Q_h^n)^2)
        &= (Q_h^{n+1})^2:I - (Q_h^n)^2:I \\
        &= (Q_h^{n+1} + Q_h^n):(Q_h^{n+1} - Q_h^n) \\
        &= 2Q_h^{n + \half}:(Q_h^{n+1}-Q_h^n),
    \end{myeq*}
    and
    \begin{myeq*}
        \tr((Q_h^{n+1})^3) - \tr((Q_h^n)^3)
        &= (Q_h^{n+1})^3:I - (Q_h^n)^3:I \\
        &= ((Q_h^{n+1})^2 + Q_h^{n+1} Q_h^n + (Q_h^n)^2):(Q_h^{n+1} - Q_h^n) \\
        &= (2(Q_h^2)^{n+\half} + Q_h^{n+1}Q_h^n):(Q_h^{n+1} - Q_h^n),
    \end{myeq*}
    where we have used that $Q^{n+1}_h$ and $Q^n_h$ are symmetric. Next, we have
    \begin{myeq*}
        \tr((Q_h^{n+1})^2)^2 - \tr((Q_h^n)^2)^2
        &= \paren{\tr((Q_h^{n+1})^2) + \tr((Q_h^n)^2)} \paren{\tr((Q_h^{n+1})^2) - \tr((Q_h^n)^2)} \\
        &= \paren{\tr((Q_h^{n+1})^2) + \tr((Q_h^n)^2)} (Q_h^{n+1} + Q_h^n):(Q_h^{n+1} - Q_h^n) \\
        &= 4\paren{\tr(Q_h^2)}^{n+\half} Q_h^{n + \half} : (Q_h^{n+1} - Q_h^n).
    \end{myeq*}
    Thus,
    \begin{myeq*}
        & \mathcal{F}_6(Q_h^{n+1}) - \mathcal{F}_6(Q_h^n) \\
        &= \myang{
            a(2Q_h^{n + \half})
            - \frac{2b}3 (2(Q_h^2)^{n+\half} + Q_h^{n+1}Q_h^n)
            + \frac c 2 (4\paren{\tr(Q_h^2)}^{n+\half} Q_h^{n + \half})
            , Q_h^{n+1} - Q_h^n
            }
            \\
        &= H_6^{n + \half}(Q_h^{n+1} - Q_h^n).
    \end{myeq*}
    We combine the above expressions and use Scheme~\ref{fullscheme} to complete the proof.
    \begin{myeq*}
        &\frac{\mathcal{F}(Q_h^{n+1}) - \mathcal{F}(Q_h^n)} {\Delta t} = \frac 1{\Delta t} H^{n+\half}(Q_h^{n+1}-Q_h^n)= H^{n + \half} \paren{\frac{Q_h^{n+1} - Q_h^n}{\Delta t}} 
        = -\frac 1 M \myVert{\frac{Q_h^{n+1} - Q_h^n}{\Delta t}}^2
    \end{myeq*}
\end{proof}

As a preliminary for the following theorem on the solvability of Scheme~\ref{fullscheme}, we define the following functionals in $\paren{\RR^{3 \times 3} \to \RR}^2 \to \RR$ for all $n$.
\begin{myeq}
    H_{1a}^{n + \half}(X,\varphi) &= \frac{L_1}4 \myang{S_1(X) \,\div X + S_1(Q_h^n) \,\div Q_h^n, (S_1(X) + S_1(Q_h^n)) \,\div \varphi} \\
    H_{1b}^{n + \half}(X,\varphi) &= \frac{L_1}4 \myang{S_1(X) \,\div X + S_1(Q_h^n) \,\div Q_h^n, \varphi (\div X + \div Q_h^n)} \\
    H_{2a}^{n + \half}(X,\varphi) &= \frac{L_2}4 \myang{S_1(X) \,\curl X + S_1(Q_h^n) \,\curl Q_h^n, (S_1(X) + S_1(Q_h^n)) \,\curl \varphi} \\
    H_{2b}^{n + \half}(X,\varphi) &= \frac{L_2}4 \myang{S_1(X) \,\curl X + S_1(Q_h^n) \,\curl Q_h^n, \varphi (\div X + \curl Q_h^n)} \\
    H_{3a}^{n + \half}(X,\varphi) &= \frac{L_3}4 \myang{S_2(X) \,\div X + S_2(Q_h^n) \,\div Q_h^n, (S_2(X) + S_2(Q_h^n)) \,\div \varphi} \\
    H_{3b}^{n + \half}(X,\varphi) &= \frac{L_3}4 \myang{S_2(X) \,\div X + S_2(Q_h^n) \,\div Q_h^n, -\varphi (\div X + \div Q_h^n)} \\
    H_{4a}^{n + \half}(X,\varphi) &= \frac{L_4}4 \myang{S_2(X) \,\curl X + S_2(Q_h^n) \,\curl Q_h^n, (S_2(X) + S_2(Q_h^n)) \,\curl \varphi} \\
    H_{4b}^{n + \half}(X,\varphi) &= \frac{L_4}4 \myang{S_2(X) \,\curl X + S_2(Q_h^n) \,\curl Q_h^n, -\varphi (\div X + \curl Q_h^n)} \\
    H_{5a}^{n + \half}(X,\varphi) &= \frac{L_5}4 \myang{\paren{|X|^2 + |Q_h^n|^2} (\nabla X + \nabla Q_h^n), \nabla \varphi} \\
    H_{5b}^{n + \half}(X,\varphi) &= \frac{L_5}4 \myang{\paren{|\nabla X|^2 + |\nabla Q_h^n|^2) (X + Q_h^n}, \varphi} \\
    H_{6a}^{n + \half}(X,\varphi) &= a \myang{X + Q_h^n,\varphi} \\
    H_{6b}^{n + \half}(X,\varphi) &= -\frac {2b}3 \myang{X^2 + (Q_h^n)^2 + XQ_h^n, \varphi} \\
    H_{6c}^{n + \half}(X,\varphi) &= \frac c2 \myang{\paren{\tr X^2 + \tr(Q_h^n)^2}(X + Q_h^n), \varphi} \\
    H^{n + \half}(X,\varphi) &= H_{1a}^{n + \half}(X,\varphi) + H_{1b}^{n + \half}(X,\varphi) + H_{2a}^{n + \half}(X,\varphi) + H_{2b}^{n + \half}(X,\varphi) + H_{3a}^{n + \half}(X,\varphi) \\&+ H_{3b}^{n + \half}(X,\varphi) + H_{4a}^{n + \half}(X,\varphi) + H_{4b}^{n + \half}(X,\varphi) + H_{5a}^{n + \half}(X,\varphi) + H_{5b}^{n + \half}(X,\varphi) \\&+ H_{6a}^{n + \half}(X,\varphi) + H_{6b}^{n + \half}(X,\varphi) + H_{6c}^{n + \half}(X,\varphi)
\end{myeq}
Now we can write the scheme as follows:
\begin{myeq}
    \myang{\frac{Q_h^{n+1} - Q_h^n}{\Delta t}, \varphi} &= -M H^{n + \half}(Q_h^{n+1}, \varphi) \quad \forall \varphi \in \mathcal T_h^0.
\end{myeq}

We write $\widetilde{Q}_h^{n+1} = Q_h^{n+1}-G_h$ and $\widetilde{Q}_h^n =Q_h^n -G_h$ which have both homogeneous Dirichlet boundary conditions and are therefore in $\mathcal{T}_h^0$ and solve
\begin{equation}\label{eq:Qtilde}
    \myang{\frac{\widetilde{Q}_h^{n+1} - \widetilde{Q}_h^n}{\Delta t}, \varphi} = -M H^{n + \half}(\widetilde{Q}_h^{n+1}+G_h, \varphi) \quad \forall \varphi \in \mathcal T_h^0.
\end{equation}
For all $n$, we consider the mapping $X\mapsto \mathcal L^{n+1}(X)$, where $\mathcal L^{n+1}(X)$ is the unique element of $\mathcal T_h^0$ such that
\begin{myeq}\label{Lgdef}
    \myang{\mathcal L^{n+1}(X), \varphi} = \myang{\widetilde{Q}_h^n,\varphi} - M \Delta t H^{n+\half}(X+G_h,\varphi) \quad \forall \varphi \in \mathcal T_h^0.
\end{myeq}
This element is uniquely defined because the mass matrix $\langle \varphi_i,\varphi_j\rangle_{ij}$ is invertible. Clearly, if this mapping has a fixed point, it is $\widetilde{Q}_h^{n+1}$ which solves~\eqref{eq:Qtilde} and then $Q^{n+1}_h$ defined as $Q^{n+1}_h=\widetilde{Q}_h^{n+1}+G_h$ corresponds to a solution of Scheme~\ref{fullscheme}. We prove this in the next theorem using the Leray-Schauder fixed point theorem~\cite[page 540]{Evans2010}. Using a contraction mapping argument, one can show that the fixed point is unique, as long as $\Delta t \leq C h^{d+2}$ for some constant $C$ independent of the meshsize $h$ and the time step size $\Delta t$.
\begin{theorem}\label{uniquesol}
    Scheme \ref{fullscheme} admits a solution $Q_h^{n+1}$ for all $n$. If $\Delta t \leq C h^{d+2}$ for some constant $C$ independent of $h$ and $\Delta t$, then the solution is unique.
\end{theorem}
\begin{proof}
We use the Leray-Schauder theorem to show existence of a fixed point $\widetilde{Q}_h^{n+1}$ of the mapping $\mathcal{L}^{n+1}:\mathcal T_h^0\to \mathcal T_h^0$ defined in~\eqref{Lgdef} above.  Then $Q^{n+1}_h:= \widetilde{Q}_h^{n+1}+G_h$ solves Scheme~\ref{fullscheme}. To apply the Leray-Schauder theorem, we need to show that the mapping is well-posed, compact, continuous,  as a mapping from $\mathcal{T}_h^0$ to $\mathcal{T}_h^0$, and that the set
		\begin{equation*}
			K:=\{Q\in \mathcal{T}_h^0\, |\, Q = \lambda \mathcal{L}^{n+1}(Q)\quad \text{for some }\ 0\leq \lambda \leq 1\}
		\end{equation*}
		is bounded~\cite[page 540]{Evans2010} (we also need $0\in \mathcal{T}_h^0$, which is clearly true). We omit the dependence of $\mathcal{L}$ on $n+1$ in the following proof, as it is not relevant.
		We assume without loss of generality that $\mathcal{F}$ is nonnegative ($\mathcal{F}$ is bounded from below, and so if it is negative, one can add a sufficiently large constant to make it positive and not change the dynamics or minimizers of the equation.).
		We start with the wellposedness.
		
		\medskip
		
	{\bf Well-posedness of the mapping $\mathcal{L}$}: This follows from above. The mass matrix $M_{ij}:=\langle \varphi_i,\varphi_j\rangle_{ij}$ is invertible and hence the solution of the mapping is uniquely defined.
		
		\medskip
		
		{\bf Continuity: } Let $\mathcal{L}(X) = \sum_{j=1}^{|\mathcal{B}_h|}\hat{X}_i \varphi_i$ and $\mathcal{L}(Y) = \sum_{j=1}^{|\mathcal{B}_h|}\hat{Y}_i \varphi_i$ two solutions of the mapping. Then we have
		\begin{equation*}
			\norm{\mathcal{L}(X)-\mathcal{L}(Y)}_{L^2}^2 =  -M\Delta t( H^{n+\frac{1}{2}}(X+G_h,\mathcal{L}(X)-\mathcal{L}(Y) ) - H^{n+\frac{1}{2}}(Y+G_h,\mathcal{L}(X)-\mathcal{L}(Y) )).
		\end{equation*}
	We show that the terms $H^{n+\frac12}_{5a}$ and $H^{n+\frac12}_{5b}$ terms are continuous, the terms $H^{n+\frac12}_{ia}$ and $H^{n+\frac12}_{ib}$, $i=1,\dots, 4$ are similar and the $H^{n+\frac12}_{6}$-terms are polynomial in $X$ and therefore continuous. So we have by the triangle inequality,
	\begin{align*}
	&	\left|H^{n+\frac12}_{5a}(X+G_h,\mathcal{L}(X)-\mathcal{L}(Y) ) - H^{n+\frac12}_{5a}(Y+G_h,\mathcal{L}(X)-\mathcal{L}(Y) )\right|\\
	&\quad  = \frac{L_5}{4}\left|\langle  (|X+G_h|^2 +|Q^n_h|^2) (\nabla X +\nabla G_h+\nabla Q_h^n)- (|Y+G_h|^2 +|Q^n_h|^2) (\nabla Y +\nabla G_h+\nabla Q_h^n), \nabla(\mathcal{L}(X)-\mathcal{L}(Y) )  \rangle \right|\\
		&\quad  \leq \frac{L_5}{4}\left|\langle (|X+G_h|^2-|Y+G_h|^2)(\nabla X +\nabla G_h +\nabla Q_h^n),\nabla(\mathcal{L}(X)-\mathcal{L}(Y)  ) \rangle \right|\\
		&\quad +\frac{L_5}{4}\left|\langle (|Y+G_h|^2+|Q_h^n|^2)(\nabla X -\nabla Y),\nabla(\mathcal{L}(X)-\mathcal{L}(Y))   \rangle \right|.
        \end{align*}
        Using H\"older's inequality for the right hand side, we obtain,
        \begin{align*}
        &	\left|H^{n+\frac12}_{5a}(X+G_h,\mathcal{L}(X)-\mathcal{L}(Y) ) - H^{n+\frac12}_{5a}(Y+G_h,\mathcal{L}(X)-\mathcal{L}(Y) )\right|\\
		&\quad \leq C \norm{|X|+|Y|+2|G_h|}_{L^\infty}\norm{X-Y}_{L^2}\norm{\nabla X +\nabla G_h +\nabla Q_h^n}_{L^\infty}\norm{\nabla(\mathcal{L}(X)-\mathcal{L}(Y))}_{L^2}\\
		&\quad +C\left(\norm{Y+G_h}_{L^\infty}^2+\norm{Q_h^n}^2_{L^\infty}\right)\norm{\nabla(X-Y)}_{L^2}\norm{\nabla(\mathcal{L}(X)-\mathcal{L}(Y)) }_{L^2}.
        \end{align*}
Next, we use inverse inequalities, to bound the last terms as follows:
        \begin{align*}
    &	\left|H^{n+\frac12}_{5a}(X+G_h,\mathcal{L}(X)-\mathcal{L}(Y) ) - H^{n+\frac12}_{5a}(Y+G_h,\mathcal{L}(X)-\mathcal{L}(Y) )\right|\\
		&\quad \leq C h^{-d-2}(\norm{X}_{L^2}+\norm{Y}_{L^2}+C)\norm{X-Y}_{L^2}\left(\norm{ X }_{L^2}+C\right)\norm{\mathcal{L}(X)-\mathcal{L}(Y)}_{L^2}\\
		&\quad +Ch^{-d-2}\left(C+\norm{Y}_{L^2}^2\right)\norm{X-Y}_{L^2}\norm{\mathcal{L}(X)-\mathcal{L}(Y) }_{L^2}\\
		& \quad \leq C h^{-d-2} \norm{X-Y}_{L^2}\norm{\mathcal{L}(X)-\mathcal{L}(Y) }_{L^2}.
	\end{align*} 
	For the $H^{n+\frac12}_{5b}$ term we have, using triangle inequality,
		\begin{align*}
		&	\left|H^{n+\frac12}_{5b}(X+G_h,\mathcal{L}(X)-\mathcal{L}(Y) ) - H^{n+\frac12}_{5b}(Y+G_h,\mathcal{L}(X)-\mathcal{L}(Y) )\right|\\
		&\quad  = \frac{L_5}{4}\left|\langle  (|\nabla(X+G_h)|^2 +|\nabla Q^n_h|^2) (X +G_h+ Q_h^n)- (|\nabla (Y+G_h)|^2 +|\nabla Q^n_h|^2) (Y + G_h+ Q_h^n), \mathcal{L}(X)-\mathcal{L}(Y)   \rangle \right|\\
		&\quad  \leq \frac{L_5}{4}\left|\langle (|\nabla (X+G_h)|^2-|\nabla(Y+G_h)|^2)( X + G_h + Q_h^n),\mathcal{L}(X)-\mathcal{L}(Y)   \rangle \right|\\
		&\quad +\frac{L_5}{4}\left|\langle (|\nabla(Y+g_h)|^2+|\nabla Q_h^n|^2)( X - Y),\mathcal{L}(X)-\mathcal{L}(Y)   \rangle \right|.
	\end{align*}
	The last expression can be bounded as follows, using H\"older's inequality,
	\begin{align*}
		&	\left|H^{n+\frac12}_{5b}(X+G_h,\mathcal{L}(X)-\mathcal{L}(Y) ) - H^{n+\frac12}_{5b}(Y+G_h,\mathcal{L}(X)-\mathcal{L}(Y) )\right|\\
		&\quad \leq C \norm{|\nabla X|+|\nabla Y|+2|\nabla G_h|}_{L^\infty}\norm{\nabla(X-Y)}_{L^2}\norm{ X + G_h + Q_h^n}_{L^\infty}\norm{\mathcal{L}(X)-\mathcal{L}(Y)}_{L^2}\\
		&\quad +C\left(\norm{
			\nabla (Y+G_h) }_{L^\infty}^2 +\norm{\nabla Q_h}^2_{L^\infty}\right)\norm{X-Y}_{L^2}\norm{\mathcal{L}(X)-\mathcal{L}(Y) }_{L^2}.
            \end{align*}
Again, we use inverse inequalities to bound the last terms on the right hand side, we have 
            \begin{align*}
                &	\left|H^{n+\frac12}_{5b}(X+G_h,\mathcal{L}(X)-\mathcal{L}(Y) ) - H^{n+\frac12}_{5b}(Y+G_h,\mathcal{L}(X)-\mathcal{L}(Y) )\right|\\
                &\quad \leq C h^{-d-2}\left(\norm{X}_{L^2}+\norm{Y}_{L^2}+2\norm{G_h}_{L^2}\right)\norm{X-Y}_{L^2}\norm{ X + G_h + Q_h^n}_{L^2}\norm{\mathcal{L}(X)-\mathcal{L}(Y)}_{L^2}\\
		&\quad + C h^{-d-2}\left(\norm{
			Y+G_h }_{L^2}^2 +\norm{Q_h^n}^2_{L^2}\right)\norm{X-Y}_{L^2}\norm{\mathcal{L}(X)-\mathcal{L}(Y) }_{L^2}\\
			&\quad \leq C h^{-d-2}\norm{X-Y}_{L^2}\norm{\mathcal{L}(X)-\mathcal{L}(Y) }_{L^2}.
	\end{align*}
	Using the same strategy for the other terms, we obtain
    \begin{equation}\label{eq:continuity}
        \norm{\mathcal{L}(X)-\mathcal{L}(Y)}_{L^2}^2\leq C\Delta t h^{-d-2}\norm{X-Y}_{L^2}\norm{\mathcal{L}(X)-\mathcal{L}(Y) }_{L^2}.
    \end{equation}
    Dividing both sides by $ \norm{\mathcal{L}(X)-\mathcal{L}(Y)}_{L^2}$ and sending $X\to Y$, we obtain the continuity of $\mathcal{L}$. Note that this estimate~\eqref{eq:continuity} implies that if $\Delta t \leq c h^{d+2}$ for a constant $c$ sufficiently small, then
    \begin{equation*}
    	\norm{\mathcal{L}(X)-\mathcal{L}(Y)}_{L^2}< \norm{X-Y}_{L^2},
    	\end{equation*} 
    and hence $\mathcal{L}$ is a contraction if $\Delta t \leq c h^{d+2}$.
	
	\medskip 
	
		{\bf Compactness:} This follows because the space $\mathcal{T}_h^0$ is finite dimensional.
		
		\medskip
		
		{\bf Boundedness of $K$:} It remains to show that the set $K$ is uniformly bounded. We consider
		\begin{equation*}
				\myang{X, \varphi} = \lambda\myang{\widetilde{Q}_h^n,\varphi} - \lambda M \Delta t H^{n+\half}(X+G_h,\varphi) \quad \forall \varphi \in \mathcal T_h^0,		
		\end{equation*}
		and use $X-\widetilde{Q}_h^n$ as a test function. Using Theorem~\ref{discrete_dissipation} and rearranging terms, we obtain
		\begin{equation*}
			\myang{X-\lambda \widetilde{Q}_h^n, X-\widetilde{Q}_h^n} = - \lambda M \Delta t H^{n+\half}(X+G_h,X-\widetilde{Q}_h^n) =-\lambda M \Delta t (\mathcal{F}(X+G_h)-\mathcal{F}(Q_h^n)).
		\end{equation*}
		Thus
		\begin{equation*}
			\norm{X}_{L^2}^2  +\lambda \norm{\widetilde{Q}_h^n}_{L^2}^2 +\lambda M \Delta t \mathcal{F}(X+G_h) =(1+\lambda)\langle \widetilde{Q}_h^n, X\rangle +\lambda M \Delta t \mathcal{F}(Q_h^n).
		\end{equation*}
		Using the Cauchy-Schwarz inequality and bounding $\lambda$ from above and below, we have
			\begin{equation*}
			\norm{X}_{L^2}^2  \leq 2 \norm{\widetilde{Q}_h^n}_{L^2}\norm{X}_{L^2}+ M \Delta t \mathcal{F}(Q_h^n)\leq 2\norm{\widetilde{Q}_h^n}_{L^2}^2 +\frac12\norm{X}_{L^2}^2+ M \Delta t \mathcal{F}(Q_h^n) .
		\end{equation*}
		Subtracting $\frac12 \norm{X}_{L^2}^2$ on both sides, we obtain the boundedness of $X$ on $L^2$ in terms of $Q^n_h$.
		Thus the conditions of the Leray-Schauder fixed point theorem are satisfied and the mapping $\mathcal{L}$ has a fixed point. If $\Delta t\leq c h^{d+2}$, then $\mathcal{L}$ is a contraction, as noted above, and therefore the fixed point is unique by the Banach fixed point theorem. Hence the scheme is well-posed.
		
	\end{proof}

\section{\texorpdfstring{$\Gamma$-convergence}{Gamma-convergence} of the discrete energies}\label{sec:gamma}

The goal of this section is to rigorously prove the $\Gamma$-convergence of the discrete energies
\begin{equation}
    \begin{split}
        \mathcal{F}_h(Q_h)&= \int_\Omega \bigg(
\frac {L_1}2 \myvert{\paren{\frac{s_0}3 I + Q_h} \div Q_h}^2 +
\frac {L_2}2 \myvert{\paren{\frac{2s_0}3 I - Q_h} \div Q_h}^2 + \frac {L_3}2 \myvert{\paren{\frac{s_0}3 I + Q_h} \curl Q_h}^2 \\ 
&\qquad  +
\frac {L_4}2 \myvert{\paren{\frac{2s_0}3 I - Q_h} \curl Q_h}^2+ \frac{L_5}{2}\left|Q_h\right|^2\left|\nabla Q_h\right|^2  + W(Q_h) \bigg) \, dx,\quad \text{for }\, Q_h\in \mathcal{T}_h^G.
    \end{split}
\end{equation}
This implies that if $Q^*_h$ is a sequence of global minimizers of the discrete energies $\mathcal{F}_h$, then any cluster point of the sequence is a global minimizer of $\mathcal{F}$. Discrete local minimizers can be obtained through the numerical scheme for the gradient flow~\eqref{eq:gradflow}, outlined in the previous section. In general, we do not know whether these discrete local minimizers are global and unique. However, as long as $Q^*$ is an isolated local minimizer of $\mathcal{F}$, then~\cite[Thm. 5.1]{Braides2014} implies that there exists a sequence $\{Q_h^*\}$, with $Q_h^*$ a local minimizer of $\mathcal{F}_h$ for $h>0$ sufficiently small.

To begin, we define the \textit{elastic energy density} of $Q$ to be
\begin{equation*}
    \sigma(Q) = L_1|S_1 \,\div Q|^2 + L_2|S_1 \,\curl Q|^2 + L_3|S_2 \,\div Q|^2 + L_4|S_2 \,\curl Q|^2 + L_5 |\nabla Q|^2 |Q|^2,
\end{equation*}
where we will assume that  $L_1,L_2,L_3,L_4,L_5>0$, and so $\mathcal{F}(Q) = \int_\Omega (\sigma(Q) + W(Q)) \, dx$. We bound the energy density as follows.

\begin{lemma}\label{nablacaplemma}
    There exist constants $C_1$ and $C_2$ such that for any $Q$, it holds that $\sigma(Q) \leq C_1 \myvert{\nabla Q}^2 + C_2 |Q|^2 |\nabla Q|^2$.
\end{lemma}
\begin{proof} 
We consider the first term $L_1|S_1\div Q|^2$, the bounds for the other terms are similar.
	We rely on Lemmas \ref{divgrad3} and \ref{curlgrad2}.
    \begin{myeq*}
        \myvert{\paren{\frac{s_0}3 I + Q} \div Q}^2
        &\leq \paren{\frac{s_0^2}9 + 2\frac{s_0}3|Q| + |Q|^2}|\div Q|^2 \\
        &\leq \paren{\frac{s_0^2}9 + 2\frac{s_0}3\paren{\frac 14 + |Q|^2} + |Q|^2}3 |\nabla Q|^2 \\
        &= 3\paren{\frac{s_0^2}9 + \frac{s_0}6} \myvert{\nabla Q}^2 + 3 \paren{\frac{2s_0}3 + 1} |Q|^2 |\nabla Q|^2.
    \end{myeq*}
\end{proof}
Note that this lemma holds for any finite $L_i$, $i=1,\dots, 5$, including when $L_i=0$ for some $i\in \{1,\dots, 5\}$.
The converse of Lemma~\ref{nablacaplemma}, i.e., coercivity of the energy, holds true as well, as long as $L_i>0$, $i=1,\dots,5$:
\begin{lemma}
    \label{lem:coercivity}
    There exists constants $C_1$ and $C_2$ depending on the coefficients $L_i$ and the boundary data $g$ such that for any $Q\in H^1(\Omega;\mathcal{S})$,
    \begin{equation*}
        \|Q\|_{H^1}^2+\int_{\Omega}|Q|^2|\nabla Q|^2 dx \leq C_1\int_{\Omega}\sigma(Q) dx + C_2.
    \end{equation*}
\end{lemma}
\begin{proof}
    \cite[Proposition 3.1]{GOLOVATY_2020} yields 
    \begin{equation*}
        \|Q\|_{H^1}^2\leq C_1\int_{\Omega}\sigma(Q) dx + C_2.
    \end{equation*}
   
    The integral $\int_{\Omega} |Q|^2 |\nabla Q|^2 dx$ is also bounded due to the contribution of the $L_5$ term in the elastic energy.
\end{proof}
Note that this lemma requires $L_5>0$ in order to upper bound the term $\int_{\Omega}|Q|^2|\nabla Q|^2 dx$. If $L_5=0$, it is unclear whether this lower bound still holds. As we shall see later, this lemma is crucial for constructing a recovery sequence.
Lemmas~\ref{nablacaplemma} and~\ref{lem:coercivity} imply the bounds
\begin{equation*}
    C_1\left(\|Q_h\|_{H^1}^2+\int_{\Omega}|Q_h|^2|\nabla Q_h|^2 dx \right) \leq \mathcal{F}_h(Q_h)\leq C_2 \left(\|Q_h\|_{H^1}^2+\int_{\Omega}|Q_h|^2|\nabla Q_h|^2 dx \right)
\end{equation*}
with $C_1, C_2>0$ independent of $h$ for the discrete energies $\mathcal{F}_h$.
To prove $\Gamma$-convergence of the discrete energies, we need to show weak lower semincontinuity (the ``liminf property'') and the existence of a recovery sequence (the ``limsup property,'' or consistency). To establish the existence of a recovery sequence, we need a few preliminary lemmas that we will prove next. First, we define the \textit{truncation} of $Q:\Omega \to \mathcal S$ by $R > 0$ as
\begin{myeq}
    T_R(Q) = \begin{cases}
        Q & \myvert{Q} \leq R \\
        \frac{R}{\myvert Q} Q & \myvert Q > R
    \end{cases}.
\end{myeq}
\begin{lemma}\label{truncgradlemma}
    For any $Q$ and $R$, it holds that $|\nabla T_R(Q)| < 2|\nabla Q|$. 
\end{lemma}
\begin{proof}
    If $|Q| \leq R$, then $T_R(Q) = Q$ and the lemma is trivial. If $|Q| > R$, then
    \begin{myeq*}
        \myvert{\nabla T_R(Q)}
        &= R \sqrt{\sum_{i,j,m=1}^3 \paren{\partial_m \paren{\frac{Q_{ij}}{\sqrt{\sum_{k,\ell=1}^3 Q_{k\ell}^2}}}}^2} \\
        &= R \sqrt{\sum_{i,j,m=1}^3 \paren{ \frac{Q_{ij,m}}{|Q|} - \frac{Q_{ij} \sum_{k,\ell=1}^3 Q_{k\ell,m} Q_{k\ell}}{|Q|^3} }^2 } \\
        &\leq \frac R{|Q|}|\nabla Q| + \frac R{|Q|^3} \sqrt{\paren{\sum_{i,j=1}^3 Q_{ij}^2} \paren{\sum_{k,\ell,m=1}^3 Q_{k\ell,m}^2} \paren{\sum_{k,\ell=1}^3 Q_{k\ell}^2} } \\
        &= 2 \frac R{|Q|}|\nabla Q|.
    \end{myeq*}
\end{proof}
%
Next, we show that the energies of the truncated $Q$-tensors approximate the standard truncated energies well for large $R>0$.
\begin{lemma}\label{gamma1}
    For every $Q \in H^1(\Omega; \mathcal S)$, we have
    \begin{myeq*}
        T_R(Q) \to Q \text{ in } H^1(\Omega; \mathcal S) \text{ and } \mathcal F(T_R(Q)) \to \mathcal F(Q) \quad \text{as }R \to \infty.
    \end{myeq*}
\end{lemma}
\begin{proof}
    $T_R(Q)$ converges to $Q$ pointwise almost everywhere, $|T_R(Q)|^2 \leq |Q|^2$ by definition, and $|\nabla T_R(Q)|^2 \leq 4|\nabla Q|^2$ by Lemma \ref{truncgradlemma}, so $T_R(Q)$ converges strongly to $Q$ in $H^1$ by the Dominated Convergence Theorem. If $\int_\Omega |Q|^2 |\nabla Q|^2 \, dx = \infty$, then $\int_\Omega |T_R(Q)|^2 |\nabla T_R(Q)|^2 \, dx \to \infty$ by Fatou's lemma, and thus $\mathcal F(T_R(Q)) \to \mathcal F(Q) = \infty$. Otherwise if $\mathcal{F}(Q)<\infty$, we have by Lemmas~\ref{nablacaplemma} and~\ref{truncgradlemma}
    \begin{myeq*}
        \sigma(T_R(Q)) &\leq C_1 |\nabla T_R(Q)|^2 + C_2 |T_R(Q)|^2 |\nabla T_R(Q)|^2 \leq 4C_1|\nabla Q|^2 + 4C_2|Q|^2 |\nabla Q|^2\leq C\sigma(Q), \\
        W(T_R(Q)) &\leq W(Q), 
    \end{myeq*}
    using Lemma~\ref{lem:coercivity} for the first term,
    so the integrand of $\mathcal F(T_R(Q))$ is also dominated by an integrable quantity, and we can use the Dominated Convergence Theorem again.
\end{proof}

\begin{lemma}\label{gamma2}
    If $\Omega$ is bounded with Lipschitz regular boundary, then for every $Q \in H^1(\Omega; \mathcal S) \cap L^{\infty}(\Omega; \mathcal S)$, there exists a sequence $\{Q^m\} \in C^\infty(\overline \Omega)$ such that
    \begin{myeq}
        Q^m \to Q \text{ in } H^1(\Omega; \mathcal S) \text{ and } \mathcal F(Q^m) \to \mathcal F(Q) \quad \text{as }m \to \infty.
    \end{myeq}
\end{lemma}

\begin{proof}
    Let $Q^m$ be the mollification of $Q$ as defined in \cite[Chapter 5.3]{Evans2010} (see Theorem 4.3 in~\cite{Evans2015} for the approximation near the boundary); then $\|Q^m\|_{L^\infty} \leq \|Q\|_{L^\infty}$, and $Q^m \to Q$ in $H^1(\Omega; \mathcal S)$ as $m \to \infty$. By considering a subsequence, we can take $Q^m$ to converge pointwise almost everywhere without loss of generality. Then we bound the difference between $\mathcal{F}_1(Q^m)$ and $\mathcal{F}_1(Q)$ as follows:
    \begin{myeq*}
        &\myvert{\mathcal{F}_1(Q^m) - \mathcal{F}_1(Q)} \leq \int_\Omega \myvert{
            \myvert{\paren{\frac{s_0}3 I + Q^m} \div Q^m}^2 - \myvert{\paren{\frac{s_0}3 I + Q} \div Q}^2
        } \, dx \\
        &= \int_\Omega \myvert{
            \paren{\paren{\frac{s_0}3 I + Q^m} \div Q^m + \paren{\frac{s_0}3 I + Q} \div Q} :
            \paren{\paren{\frac{s_0}3 I + Q^m} \div Q^m - \paren{\frac{s_0}3 I + Q} \div Q}} \, dx \\
        &\leq \int_\Omega 
            \paren{\myvert{\paren{\frac{s_0}3 I + Q^m} \div Q^m} + \myvert{\paren{\frac{s_0}3 I + Q} \div Q}}
            \myvert{\paren{\frac{s_0}3 I + Q^m} \div Q^m - \paren{\frac{s_0}3 I + Q} \div Q}
          dx.
       \end{myeq*}
       Using the Cauchy-Schwarz inequality for the last expression, we obtain
       \begin{myeq*}
             &\myvert{\mathcal{F}_1(Q^m) - \mathcal{F}_1(Q)} 	\\
        &\leq \paren{\paren{\mathcal{F}_1(Q^m)}^\half + \paren{\mathcal{F}_1(Q)}^\half} \paren{\int_\Omega \myvert{\paren{\frac{s_0}3 I + Q^m} \div Q^m - \paren{\frac{s_0}3 I + Q} \div Q}^2 \, dx}^\half .
    \end{myeq*}
   Next, we use the triangle inequality, to get
    \begin{myeq*}
    &\myvert{\mathcal{F}_1(Q^m) - \mathcal{F}_1(Q)} 	\\
        &= \paren{\paren{\mathcal{F}_1(Q^m)}^\half + \paren{\mathcal{F}_1(Q)}^\half} \paren{\int_\Omega \myvert{\paren{\frac{s_0}3 I + Q^m} \div (Q^m - Q) + (Q^m - Q) \, \div Q}^2 \, dx}^\half \\
        &\leq \paren{\paren{\mathcal{F}_1(Q^m)}^\half + \paren{\mathcal{F}_1(Q)}^\half} \left(\left(\int_\Omega \myvert{\paren{\frac{s_0}3 I + Q^m} \div (Q^m - Q)}^2 \, dx \right)^{\frac12}
            + \left(  \int_\Omega \myvert{(Q^m - Q) \, \div Q}^2 \, dx .
        \right)^{\frac12} \right)
    \end{myeq*}
    To prove that $\mathcal{F}_1(Q^m)$ is uniformly bounded by a constant, we observe that
    \begin{myeq*}
        \int_\Omega \myvert{\paren{\frac{s_0}3 I + Q^m} \div Q^m}^2 \, dx
        \leq \int_\Omega \myvert{\frac{s_0}3 I + Q^m}^2 \myvert{\div Q^m}^2 \, dx
        \leq 9 \myVert{\frac{s_0}3 I + Q^m}_{L^\infty}^2 \int_\Omega 3 \myvert{\nabla Q^m}^2 \, dx,
    \end{myeq*}
    which is bounded because $\myVert{Q^m}_{L^\infty} \leq \myVert{Q}_{L^\infty}$ and $Q^m\to Q$ in $H^1(\Omega; \mathcal S)$ and thus $\|Q^m\|_{H^1}\leq \|Q\|_{H^1}+C$ for $m$ large enough.
    For the first remaining integral, we have
    \begin{myeq*}
        \int_\Omega \myvert{\paren{\frac{s_0}3 I + Q^m} \div (Q^m - Q)}^2 \, dx &\leq 9\myVert{\frac{s_0}3 I + Q^m}_{L^\infty}^2 \int_\Omega |\div(Q^m - Q)|^2,
    \end{myeq*}
    which approaches zero as $m \to \infty$ due to the $H^1$-convergence of $Q^m$ to $Q$. For the second remaining integral, we observe that $|(Q^m - Q) \, \div Q|^2$ approaches zero almost everywhere in $\Omega$, and
    \begin{myeq*}
        |(Q^m - Q) \, \div Q|^2 \leq |Q^m - Q|^2 |\div Q|^2 \leq (|Q^m| + |Q|)^2 |\div Q|^2,
    \end{myeq*}
    which is integrable; thus by the Dominated Convergence Theorem, this integral also approaches zero as $m \to \infty$. Therefore, $\mathcal{F}_1(Q^m) \to \mathcal{F}_1(Q)$. The proof for the other energy terms is similar.
\end{proof}
We have two definitions to make before stating the main theorem of this section. First, let $\{\Omega^h\}_{h > 0}$ be a regular family of triangulations of $\Omega$ such that for every $h$, each element in $\Omega^h$ has diameter less than $h$.  Second, for $Q \in C(\overline{\Omega}; \mathcal S)$ and all $h$, let $\mathcal L_h(Q)$ be the piecewise linear Lagrange interpolation of $Q$ onto $\Omega^h$. 
Then we have
\begin{lemma}\label{gamma3}
    For every $Q \in C^\infty(\overline \Omega)$, we have
    \begin{myeq}
        \myvert{\mathcal F(\mathcal{L}_h (Q)) - \mathcal F(Q)} \to 0 \quad \text{as } h \to 0.
    \end{myeq}
\end{lemma}
\begin{proof}
    From Theorem 3.2 of \cite{Bartels2015}, we have that $\mathcal L_h(Q) \to Q$ in $H^1(\Omega; \mathcal S)$ as $h \to 0$. After this observation, we can just follow the proof of Lemma~\ref{gamma2}.
\end{proof}

\begin{theorem}\label{thm:gamma}
    Let $\Omega$ be bounded with Lipschitz regular boundary; then the sequence $\{\mathcal F_h\}$ is $\Gamma$-convergent in the weak topology of $H^1(\Omega;\mathcal S)$. That is,
    \begin{enumerate}
        \item For any $Q \in H^1(\Omega;\mathcal S)$ and for any sequence $\{Q_h\}$ in $H^1(\Omega;\mathcal S)$,
        \begin{myeq}
            Q_h \rightharpoonup Q \text{ in } H^1(\Omega; \mathcal S) \text{ implies } \liminf_{h \to 0} \mathcal F_h(Q_h) \geq \mathcal F(Q),
        \end{myeq}
         and
         \item For each $Q \in H^1(\Omega;\mathcal S)$, there exists a recovery sequence $\{Q_h\}$ in $H^1(\Omega;\mathcal S)$ satisfying
         \begin{myeq}
             Q_h \rightharpoonup Q \text{ in } H^1(\Omega; \mathcal S) \text{ and } \lim_{h \to 0} \mathcal F_h(Q_h) = \mathcal F(Q).
         \end{myeq}
    \end{enumerate}
\end{theorem}

\begin{proof}
The lower semicontinuity (liminf property) follows from Proposition 3.2 in~\cite{GOLOVATY_2020}. 
To begin construction of a recovery sequence, take $k \in \ZZ^+$. By Lemma \ref{gamma1}, there exists a sequence $\{R_k\}_k$ and $\overline{R}$ large enough such that for all $R_k\geq \overline{R}$,
\begin{myeq}
    \myVert{T_{R_k}(Q) - Q}_{H^1} < \frac 1{2k} \text{ and } \begin{cases}
        \myvert{\mathcal F(T_{R_k}(Q)) - \mathcal F(Q)} < \frac 1{3k} & \text{if } \mathcal F(Q) < \infty, \\
        \mathcal F(T_{R_k}(Q)) > k + \frac{2}{3k} & \text{if } \mathcal F(Q) = \infty.
    \end{cases}
\end{myeq}
By Lemma~\ref{gamma2}, for each $R_k$ there exists a sequence $\{(T_{R_k}(Q))^{m_j}\}_{j\in\mathbb{N}} \in C^\infty(\overline \Omega)$ and $M_k>0$ such that for all $m_j\geq M_k$
\begin{myeq}
    \myVert{(T_{R_k}(Q))^{m_j} - T_{R_k}(Q)}_{H^1} < \frac 1{2j} \text{ and } \myvert{\mathcal F((T_{R_k}(Q))^{m_j}) - \mathcal F(T_{R_k}(Q))} < \frac 1{3j}.
\end{myeq}
Then, by Lemma \ref{gamma3}, for each $R_k$ and $m_j$, there exists a sequence $\{h_\ell\}_\ell > 0$ and $\bar{h}>0$ such that for all $h_\ell< \bar{h}$
\begin{myeq}
    \myvert{\mathcal F_{h_\ell}(\mathcal{L}_{h_\ell}((T_{R_k}(Q))^{m_j})) - \mathcal F((T_{R_k}(Q))^{m_j})} < \frac 1{3\ell} \text{ and } h_\ell < \begin{cases}
        1 & \text{if } \ell = 1, \\
        \min\{\frac 1\ell, h_{\ell-1}\} & \text{if } \ell > 1.
    \end{cases}
\end{myeq}
Define $Q_{h_k} = \mathcal{L}_{h_k}((T_{R_k}(Q))^{m_k})$ for all $k$. Then by the triangle inequality, we have
\begin{myeq*}
    \lim_{k \to \infty} \myVert{Q_{h_k} - Q}_{H^1} = 0 \text{ and } \lim_{k \to \infty} \mathcal F_{h_k}(Q_{h_k}) = \mathcal F(Q),
\end{myeq*}
and therefore, for $h \in \{h_1, h_2, \ldots\}$,
\begin{myeq*}
    \lim_{h \to 0} \myVert{Q_{h} - Q}_{H^1} = 0 \text{ and } \lim_{h \to 0} \mathcal F_{h}(Q_{h}) = \mathcal F(Q).
\end{myeq*}
This proves the existence of a recovery sequence for each $Q\in H^1(\Omega;\mathcal{S})$.
\end{proof}
\begin{Remark}
	We note that in Theorem~\ref{thm:gamma}, the lower semicontinuity follows as in Proposition 3.2 in~\cite{GOLOVATY_2020}, where $L_5=0$. Therefore the condition $L_5>0$ is not necessary to derive the liminf property. However, it is needed for constructing the recovery sequence, which relies on Lemma~\ref{gamma1} that bounds the energy density of the truncated $Q$, which again relies on Lemma~\ref{lem:coercivity} that is only true if $L_5>0$ as remarked above. We have experimented with other truncations of $Q$ and the elastic energy density $\sigma(Q)$ in order to avoid the constraint $L_5>0$, but were not able to remove it. However, as explained in the introduction, the $L_5$-term can be justified in the same way as the other elastic energy terms, and as the numerical experiments in the upcoming section demonstrate, the numerical results are similar for $L_5>0$ and $L_5=0$ (in particular, the numerical scheme is stable independently of the positivity of $L_5$).
\end{Remark}
\section{Numerical results in 2D}\label{sec:res}
We now present numerical simulations in two dimensions, representing a thin nematic film; the reduction from three dimensions to two is derived in ~\cite{golovaty2019phase}. We use the parameters
\begin{myeq}\label{params}
    L_1 = 0.1, \quad L_2 = L_3 = L_4 = L_5 = 0.001, \quad a = -0.3, \quad b = -4.0, \quad c = 4.0, \quad M = 1
\end{myeq}
for all experiments unless specified otherwise. To approximate $Q_h^{n+1}$, we calculate a sequence $\{Q_h^{n+1,m}\}$ using Newton iteration as follows:
    \begin{myeq}\label{eq:newtoniteration}
        D \mathcal L^{n+1}(Q_h^{n+1,m-1} - G_h) \cdot (Q_h^{n+1,m} - Q_h^{n+1,m-1}) =  -  (\mathcal L^{n+1}(Q_h^{n+1,m-1}) - Q_h^n),
    \end{myeq}
where $\mathcal L^{n+1}$ is as defined in \eqref{Lgdef} and $D \mathcal L^{n+1}$ is its Jacobian. We set $Q_h^{n+1}$ equal to the first $Q_h^{n+1,m}$ for which
\begin{myeq}\label{eq:computeepsilon}
    \myVert{Q_h^{n+1,m} - Q_h^{n+1,m-1}} < 10^{-10}.
\end{myeq}
This scheme appears to be very robust with respect to the time step size; even when we increased $\Delta t$ by several orders of magnitude in our experiments, eventually some $m$ satisfied \eqref{eq:computeepsilon} for all $n$.

These experiments were written in Julia. The code used to run them can be found at \url{https://github.com/elafandi/quartic-q-tensors}.

\subsection{Convergence tests}
To identify the order of accuracy of the scheme with respect to $h$ and $\Delta t$, we first adapt the convergence tests from \cite{gudibanda2022convergence}. For these experiments, we take our domain to be $\Omega = [0, 2]^2$ and use initial and boundary conditions
\begin{myeq}\label{convtest_Q}
    Q_0 = \mathbf{n}_0 \mathbf{n}_0^\top - \frac{|\mathbf{n}_0|^2}2 I_2,
\end{myeq}
where
\begin{myeq}\label{convtest_n}
    \mathbf{n}_0 = \begin{pmatrix}
        x(2-x) y(2-y) \\
        \sin(\pi x) \sin(\pi y/2)
    \end{pmatrix},
\end{myeq}
and 
\begin{equation*}
	\left. Q\right|_{\partial \Omega} = G =0.
\end{equation*}
For each $h$, we construct $\Omega_h$ by taking $8/h^2$ isosceles right triangles with legs of length $h$ whose legs are parallel to the axes.

\subsubsection{Refinement in space}
For $h = 0.2$, $0.1$, $0.05$, and $0.025$, corresponding to $L = 200$, $800$, $3200$, and $12800$, we compute $Q_h$ up to $T = 0.8$ with $\Delta t = 1.0 \times 10^{-3}$. We estimate the error by comparing with a reference solution $Q_{\text{ref}}$, which is calculated with $h = 0.005$ and the same $\Delta t$. Table \ref{tab:space_refine} presents the $L^2$ errors of $Q_h$ at time $T$, the absolute error of $\mathcal F[Q_h]$ at time $T$, and the rate of convergence for both. For both errors, we observe second-order convergence with respect to element width, which implies first-order convergence with respect to the number of elements.
\begin{table}[ht]
    \centering
    \begin{tabular}{|c||c|c||c|c|}
        \hline
        $h$ & $\myVert{Q_h - Q_{\text{ref}}}$ & Order & $\myvert{\mathcal F({Q_h}) - \mathcal F({Q_{\text{ref}}})}$ & Order \\
        \hline
        $0.2$ & $4.3479 \times 10^{-2}$ & $-$ & $9.8143 \times 10^{-4}$ & $-$ \\
        $0.1$ & $1.5893 \times 10^{-2}$ & $1.4519$ & $3.2932 \times 10^{-4}$ & $1.5754$ \\
        $0.05$ & $3.5399 \times 10^{-3}$ & $2.1666$ & $7.3958 \times 10^{-5}$ & $2.1547$ \\
        $0.025$ & $8.3812 \times 10^{-4}$ & $2.0785$ & $1.7802 \times 10^{-5}$ & $2.0547$ \\
        \hline
    \end{tabular}
    \caption{Errors and convergence rates for spatial refinement in (\ref{convtest_Q}), (\ref{convtest_n})}
    \label{tab:space_refine}
\end{table}

\subsubsection{Refinement in time}
We now take $h = 2/30$, and compute $Q_h$ up to $T = 0.8$ using $200$, $400$, $800$, $1600$, and $3200$ time steps. Our reference solution $Q_{\text{ref}}$ is taken with 80000 time steps and the same $h$. Results for this experiment are presented in Table \ref{tab:time_refine}. We observe second-order convergence with respect to the length of the time step as expected.
\begin{table}[h]
    \centering
    \begin{tabular}{|c||c|c||c|c|}
        \hline
        $\Delta t$ & $\myVert{Q_h - Q_{\text{ref}}}$ & Order & $\myvert{\mathcal F({Q_h}) - \mathcal F({Q_{\text{ref}}})}$ & Order \\
        \hline
        $4 \times 10^{-3}$ & $4.2759 \times 10^{-6}$ & $-$ & $4.6999 \times 10^{-7}$ & $-$ \\
        $2 \times 10^{-3}$ & $1.0688 \times 10^{-6}$ & $2.0002$ & $1.1748 \times 10^{-7}$ & $2.0002$ \\
        $1 \times 10^{-3}$ & $2.6717 \times 10^{-7}$ & $2.0001$ & $2.9368 \times 10^{-8}$ & $2.0001$ \\
        $5 \times 10^{-4}$ & $6.6772 \times 10^{-8}$ & $2.0004$ & $7.3398 \times 10^{-9}$ & $2.0004$ \\
        $2.5 \times 10^{-4}$ & $1.6673 \times 10^{-8}$ & $2.0017$ & $1.8327 \times 10^{-9}$ & $2.0018$ \\
        \hline
    \end{tabular}
    \caption{Errors and convergence rates for time refinement in (\ref{convtest_Q}), (\ref{convtest_n})}
    \label{tab:time_refine}
\end{table}

\subsection{Tactoid and nematic bubble simulations}
Inspired by \cite{golovaty2019phase}, we now simulate phase transitions in tactoids of varying topological degrees and a nematic bubble of topological degree zero. We take $\Omega$ to be the unit circle and use a Delaunay triangulation generated by \cite{Persson2004}, with $N = 5559$, $V = 5809$, and $L = 11366$. For the time step, we take $\Delta t = 0.1.$

We define $(r(x,y), \theta(x,y))$ to be the polar coordinates of $(x,y) \in \RR^2$. As an initial condition, we then take
\begin{myeq}\label{tactoid_boundary}
    Q_0  = \sqrt{\frac{-2a}{c}} \paren{ \mathbf{n}_0(\theta) \mathbf{n}_0^\top(\theta) - \frac 12 I_2 } 
    \begin{cases}
        \mathbf{1}_{r^2 \geq 0.3} & \text{for the tactoids} \\
        \mathbf{1}_{r^2 \leq 0.3} & \text{for the nematic bubble}
    \end{cases},
\end{myeq}
where $\mathbf{n}_0 : [0,2\pi) \to \mathbb{S}^1$ depends on the experiment, and set $G=\left. Q_0\right|_{\partial\Omega}$, so as to obtain a stable isotropic tactoid surrounded by a nematic sample (or reverse). We analyze two variations of each experiment: one with the parameters as given in \eqref{params}, and one where $L_1$ is set equal to $0.001$. We present the energy decay of all eight simulations in Figure \ref{fig:tactoid_energy}, and note that energy decays monotonically as expected. In addition, we compare each experiment to a simulation of the gradient flow for the standard Landau-de Gennes model with the one-constant approximation:
\begin{equation}\label{eq:oneconstLdG}
	Q_t = M\left(L\Delta Q -2a Q+2b\left(Q^2-\frac1{2}\tr( Q^2)I\right)-2c\tr(Q^2)Q\right),
\end{equation}
 with the same initial and boundary condition, where we set $L=0.001$ or $L=0.0001$ depending on the experiment and otherwise use the parameters above in~\eqref{params}.

\subsubsection{Degree 1 tactoid}\label{sec:tactOne}
\begin{figure}
    \centering
    \begin{subfigure}{0.3\textwidth}
        \includegraphics[width=\textwidth]{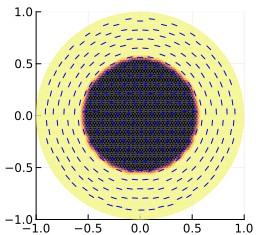}
        \caption{$t = 0.1$}
    \end{subfigure}
    \begin{subfigure}{0.3\textwidth}
        \includegraphics[width=\textwidth]{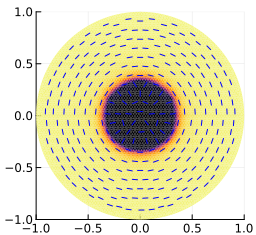}
        \caption{$t = 15$}
    \end{subfigure}
    \begin{subfigure}{0.3\textwidth}
        \includegraphics[width=\textwidth]{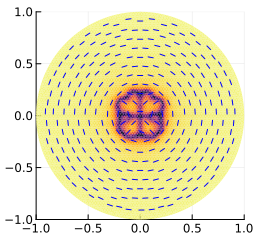}
        \caption{$t = 30$}
    \end{subfigure}
    \begin{subfigure}{0.3\textwidth}
        \includegraphics[width=\textwidth]{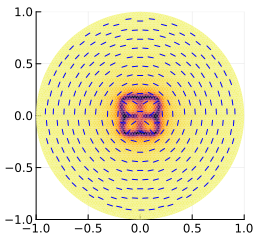}
        \caption{$t = 45$}
    \end{subfigure}
    \begin{subfigure}{0.3\textwidth}
        \includegraphics[width=\textwidth]{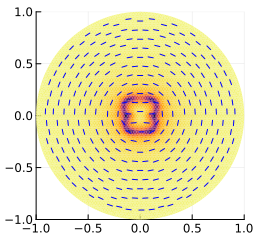}
        \caption{$t = 60$}
    \end{subfigure}
    \begin{subfigure}{0.3\textwidth}
        \includegraphics[width=\textwidth]{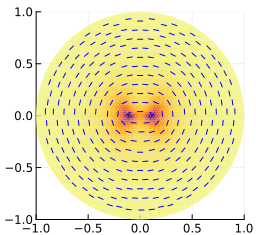}
        \caption{$t = 75$}
    \end{subfigure}
    \begin{subfigure}{0.3\textwidth}
        \includegraphics[width=\textwidth]{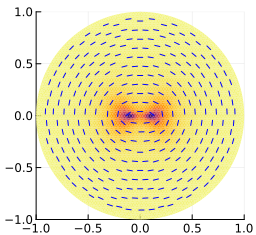}
        \caption{$t = 90$}
    \end{subfigure}
    \begin{subfigure}{0.3\textwidth}
        \includegraphics[width=\textwidth]{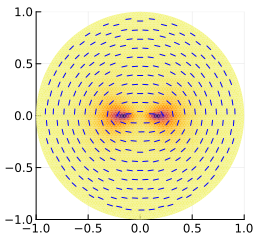}
        \caption{$t = 180$}
    \end{subfigure}
    \begin{subfigure}{0.3\textwidth}
        \includegraphics[width=\textwidth]{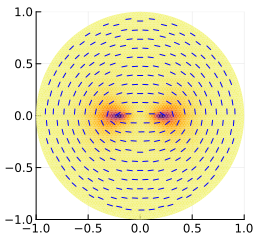}
        \caption{$t = 270$}
    \end{subfigure}
    \begin{subfigure}{0.55\textwidth}
        \includegraphics[width=\textwidth]{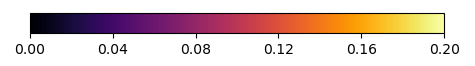}
    \end{subfigure}
    \caption{Simulated evolution of a degree 1 tactoid, in which $L_1 = 0.1$ and the initial and boundary conditions are defined by $\mathbf{n}_0 = (-\sin \theta, \cos \theta)$ in (\ref{tactoid_boundary}). Positive eigenvalues and associated eigenvectors, which correspond to the scalar order parameter and the director, are represented by the heat map and quiver plot.}
    \label{fig:tactOne_largeL1}
\end{figure}

\begin{figure}
    \centering
    \begin{subfigure}{0.3\textwidth}
        \includegraphics[width=\textwidth]{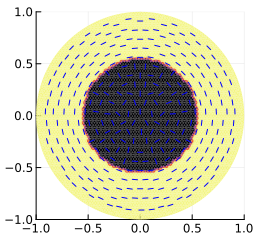}
        \caption{$t = 0.1$}
    \end{subfigure}
    \begin{subfigure}{0.3\textwidth}
        \includegraphics[width=\textwidth]{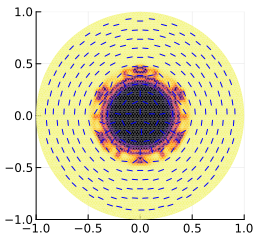}
        \caption{$t = 30$}
    \end{subfigure}
    \begin{subfigure}{0.3\textwidth}
        \includegraphics[width=\textwidth]{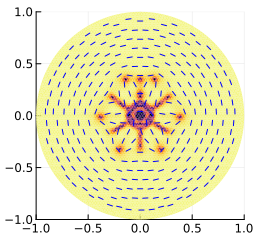}
        \caption{$t = 60$}
    \end{subfigure}
    \begin{subfigure}{0.3\textwidth}
        \includegraphics[width=\textwidth]{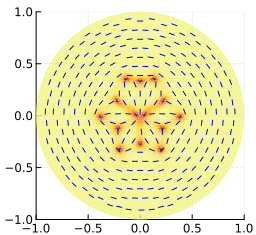}
        \caption{$t = 90$}
    \end{subfigure}
    \begin{subfigure}{0.3\textwidth}
        \includegraphics[width=\textwidth]{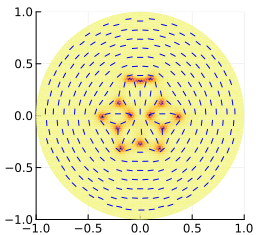}
        \caption{$t = 160$}
    \end{subfigure}
    \begin{subfigure}{0.3\textwidth}
        \includegraphics[width=\textwidth]{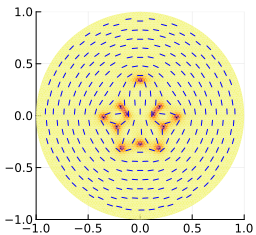}
        \caption{$t = 230$}
    \end{subfigure}
    \begin{subfigure}{0.3\textwidth}
        \includegraphics[width=\textwidth]{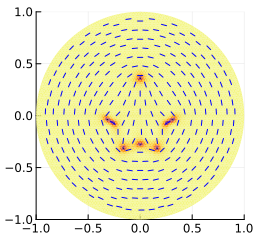}
        \caption{$t = 300$}
    \end{subfigure}
    \begin{subfigure}{0.3\textwidth}
        \includegraphics[width=\textwidth]{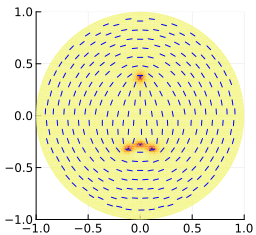}
        \caption{$t = 450$}
    \end{subfigure}
    \begin{subfigure}{0.3\textwidth}
        \includegraphics[width=\textwidth]{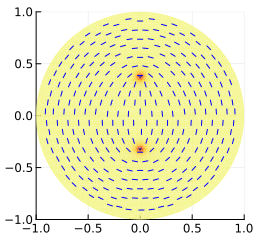}
        \caption{$t = 600$}
    \end{subfigure}
    \begin{subfigure}{0.55\textwidth}
        \includegraphics[width=\textwidth]{img/colorbar.png}
    \end{subfigure}
    \caption{Simulated evolution of a degree 1 tactoid, in which $L_1 = 0.001$ and the initial and boundary conditions are defined by $\mathbf{n}_0 = (-\sin \theta, \cos \theta)$ in (\ref{tactoid_boundary}). Positive eigenvalues and associated eigenvectors, which correspond to the scalar order parameter and the director, are represented by the heat map and quiver plot.}
    \label{fig:tactOne_smallL1}
\end{figure}

\begin{figure}
	\centering
	\begin{subfigure}{0.3\textwidth}
		\includegraphics[width=\textwidth]{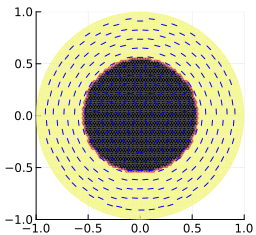}
		\caption{$t = 0.01$}
	\end{subfigure}
	\begin{subfigure}{0.3\textwidth}
		\includegraphics[width=\textwidth]{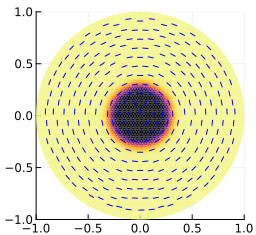}
		\caption{$t = 15$}
	\end{subfigure}
	\begin{subfigure}{0.3\textwidth}
		\includegraphics[width=\textwidth]{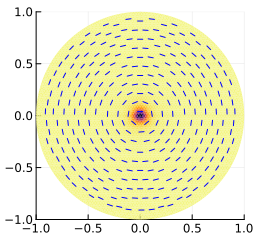}
		\caption{$t = 30$}
	\end{subfigure}
	\begin{subfigure}{0.3\textwidth}
		\includegraphics[width=\textwidth]{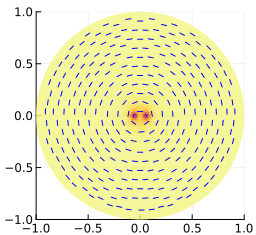}
		\caption{$t = 45$}
	\end{subfigure}
	\begin{subfigure}{0.3\textwidth}
		\includegraphics[width=\textwidth]{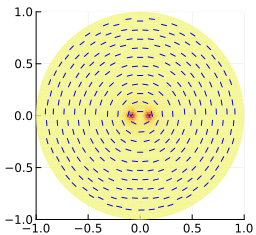}
		\caption{$t = 60$}
	\end{subfigure}
	\begin{subfigure}{0.3\textwidth}
		\includegraphics[width=\textwidth]{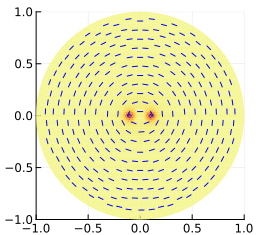}
		\caption{$t = 75$}
	\end{subfigure}
	\begin{subfigure}{0.3\textwidth}
		\includegraphics[width=\textwidth]{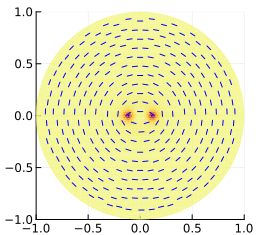}
		\caption{$t = 90$}
	\end{subfigure}
	\begin{subfigure}{0.3\textwidth}
		\includegraphics[width=\textwidth]{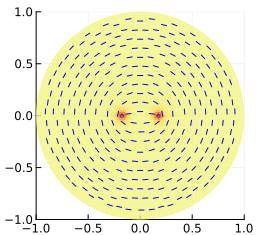}
		\caption{$t = 180$}
	\end{subfigure}
	\begin{subfigure}{0.3\textwidth}
		\includegraphics[width=\textwidth]{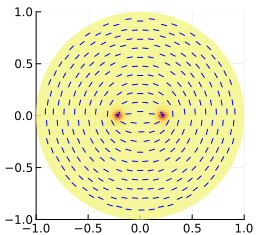}
		\caption{$t = 270$}
	\end{subfigure}
	\begin{subfigure}{0.55\textwidth}
		\includegraphics[width=\textwidth]{img/colorbar.png}
	\end{subfigure}
	\caption{Simulation of~\eqref{eq:oneconstLdG}, in which $L= 0.0001$ and the initial and boundary conditions are defined by $\mathbf{n}_0 = (-\sin \theta, \cos \theta)$ in \eqref{tactoid_boundary}. Positive eigenvalues and associated eigenvectors, which correspond to the scalar order parameter and the director, are represented by the heat map and quiver plot.}
	\label{fig:LdGdeg1}
\end{figure}

In this first trial, we take the director to be parallel to the boundary of $\Omega$ by setting
\begin{myeq}
    \mathbf{n}_0 = \begin{pmatrix} -\sin \theta \\ \cos \theta \end{pmatrix}.
\end{myeq}
This leads to a tactoid of topological degree 1: along a path which winds once clockwise around the tactoid, the director also winds once clockwise. In Figure \ref{fig:tactOne_largeL1}, we observe that when $L_1$ is large, the tactoid shrinks at an equal rate in each direction until it reaches a certain radius, at which point it slowly hollows out and eventually splits into two vortices of degree 1/2. Over a much longer time scale, these vortices drift apart from each other. Figure \ref{fig:tactOne_smallL1} shows that when $L_1$ is small, the tactoid instead decays into a constellation of vortices, all but two of which draw together and annihilate each other. In Figure~\ref{fig:LdGdeg1} we show a simulation for the same initial data for the gradient flow for the standard quadratic Landau-de Gennes model~\eqref{eq:oneconstLdG} for $L=0.0001$. In this case the spherical isotropic region shrinks in a spherically symmetric fashion until around time $t=35$ after which it splits into two 1/2 defects, which slowly move away from each other. It appears that it converges to a similar minimizer as in the case when $L_1$ is large, c.f., Figure~\ref{fig:tactMinusOne_largeL1}, however, the dynamics to reach that final state are different, the isotropic region stays spherically symmetric until the defects split up.

\begin{figure}
    \centering
    \begin{subfigure}{0.3\textwidth}
        \includegraphics[width=\textwidth]{img/rev2imgs/degree_one_L1=0.1/tactoid_deg=1_paper_30.0.png}
        \caption{$t = 30$, $L_5 = 0.001$}
    \end{subfigure}
    \begin{subfigure}{0.3\textwidth}
        \includegraphics[width=\textwidth]{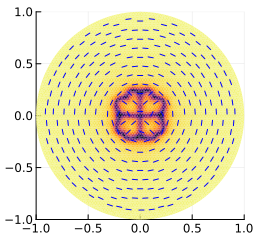}
        \caption{$t = 30$, $L_5 = 0$}
    \end{subfigure}
    \\
    \begin{subfigure}{0.3\textwidth}
        \includegraphics[width=\textwidth]{img/rev2imgs/degree_one_L1=0.1/tactoid_deg=1_paper_60.0.png}
        \caption{$t = 60$, $L_5 = 0.001$}
    \end{subfigure}
    \begin{subfigure}{0.3\textwidth}
        \includegraphics[width=\textwidth]{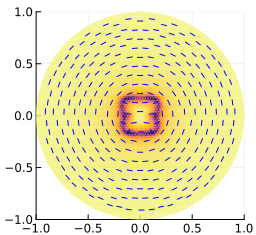}
        \caption{$t = 60$, $L_5 = 0$}
    \end{subfigure}
    \\
    \begin{subfigure}{0.3\textwidth}
        \includegraphics[width=\textwidth]{img/rev2imgs/degree_one_L1=0.1/tactoid_deg=1_paper_270.0.png}
        \caption{$t = 270$, $L_5 = 0.001$}
    \end{subfigure}
    \begin{subfigure}{0.3\textwidth}
        \includegraphics[width=\textwidth]{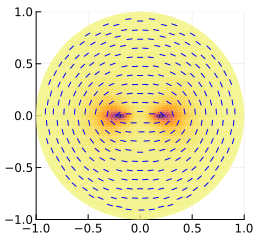}
        \caption{$t = 270$, $L_5 = 0$}
    \end{subfigure}
    \begin{subfigure}{0.55\textwidth}
        \includegraphics[width=\textwidth]{img/colorbar.png}
    \end{subfigure}
    \caption{Comparison of two simulations of a degree 1 tactoid, in which $L_1 = 0.1$ and the initial and boundary conditions are defined by $\mathbf{n}_0 = (-\sin \theta, \cos \theta)$ in (\ref{tactoid_boundary}). Positive eigenvalues and associated eigenvectors, which correspond to the scalar order parameter and the director, are represented by the heat map and quiver plot.}
    \label{fig:l5_zero}
\end{figure}

Additionally, Figure \ref{fig:l5_zero} compares the former of these simulations to one in which $L_5 = 0$. The visual differences between the two are minimal, evidencing that the $L_5$ term introduced to prove $\Gamma$-convergence does not meaningfully affect the result when sufficiently small.

\subsubsection{Degree -1 tactoid}\label{sec:tactMinusOne}
\begin{figure}
    \centering
    \begin{subfigure}{0.3\textwidth}
        \includegraphics[width=\textwidth]{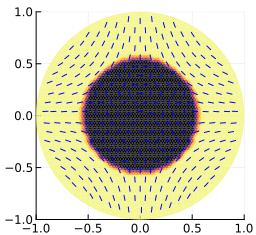}
        \caption{$t = 0.1$}
    \end{subfigure}
    \begin{subfigure}{0.3\textwidth}
        \includegraphics[width=\textwidth]{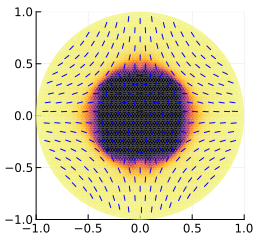}
        \caption{$t = 5$}
    \end{subfigure}
    \begin{subfigure}{0.3\textwidth}
        \includegraphics[width=\textwidth]{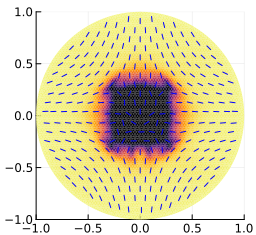}
        \caption{$t = 10$}
    \end{subfigure}
    \begin{subfigure}{0.3\textwidth}
        \includegraphics[width=\textwidth]{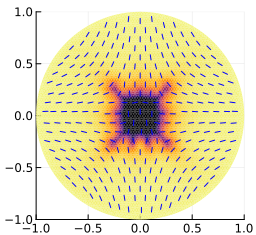}
        \caption{$t = 15$}
    \end{subfigure}
    \begin{subfigure}{0.3\textwidth}
        \includegraphics[width=\textwidth]{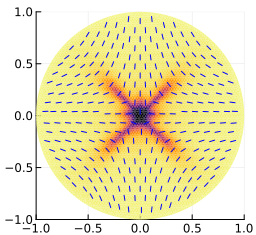}
        \caption{$t = 20$}
    \end{subfigure}
    \begin{subfigure}{0.3\textwidth}
        \includegraphics[width=\textwidth]{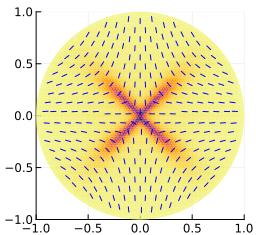}
        \caption{$t = 25$}
    \end{subfigure}
    \begin{subfigure}{0.3\textwidth}
        \includegraphics[width=\textwidth]{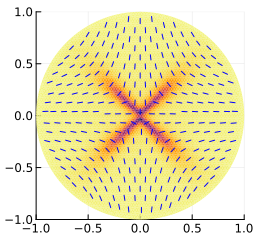}
        \caption{$t = 30$}
    \end{subfigure}
    \begin{subfigure}{0.3\textwidth}
        \includegraphics[width=\textwidth]{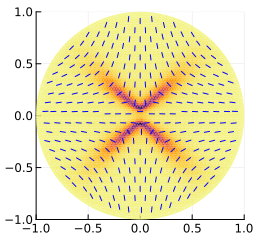}
        \caption{$t = 60$}
    \end{subfigure}
    \begin{subfigure}{0.3\textwidth}
        \includegraphics[width=\textwidth]{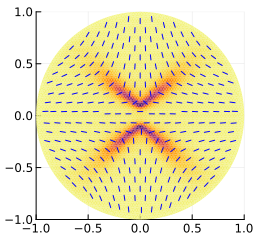}
        \caption{$t = 90$}
    \end{subfigure}
    \begin{subfigure}{0.55\textwidth}
        \includegraphics[width=\textwidth]{img/colorbar.png}
    \end{subfigure}
    \caption{Simulated evolution of a degree -1 tactoid, in which $L_1 = 0.1$ and the initial and boundary conditions are defined by $\mathbf{n}_0 = (-\cos \theta, \sin \theta)$ in (\ref{tactoid_boundary}). Positive eigenvalues and associated eigenvectors, which correspond to the scalar order parameter and the director, are represented by the heat map and quiver plot.}
    \label{fig:tactMinusOne_largeL1}
\end{figure}

\begin{figure}
    \centering
    \begin{subfigure}{0.3\textwidth}
        \includegraphics[width=\textwidth]{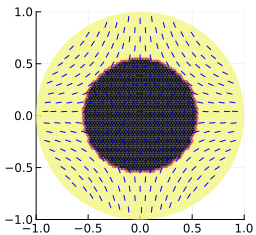}
        \caption{$t = 0.1$}
    \end{subfigure}
    \begin{subfigure}{0.3\textwidth}
        \includegraphics[width=\textwidth]{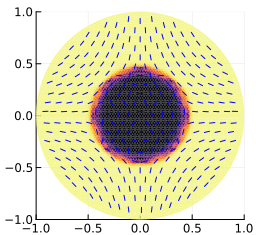}
        \caption{$t = 20$}
    \end{subfigure}
    \begin{subfigure}{0.3\textwidth}
        \includegraphics[width=\textwidth]{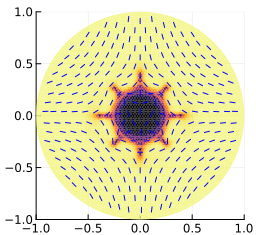}
        \caption{$t = 40$}
    \end{subfigure}
    \begin{subfigure}{0.3\textwidth}
        \includegraphics[width=\textwidth]{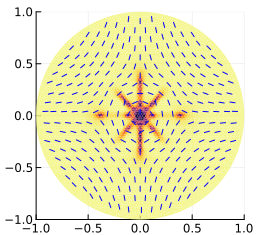}
        \caption{$t = 60$}
    \end{subfigure}
    \begin{subfigure}{0.3\textwidth}
        \includegraphics[width=\textwidth]{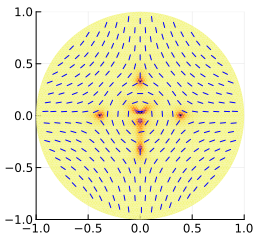}
        \caption{$t = 80$}
    \end{subfigure}
    \begin{subfigure}{0.3\textwidth}
        \includegraphics[width=\textwidth]{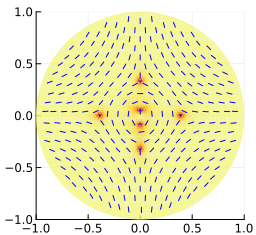}
        \caption{$t = 100$}
    \end{subfigure}
    \begin{subfigure}{0.3\textwidth}
        \includegraphics[width=\textwidth]{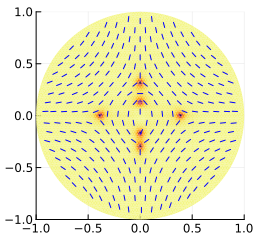}
        \caption{$t = 200$}
    \end{subfigure}
    \begin{subfigure}{0.3\textwidth}
        \includegraphics[width=\textwidth]{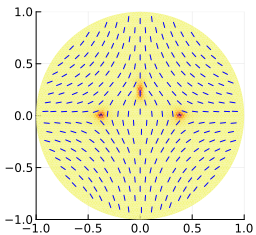}
        \caption{$t = 300$}
    \end{subfigure}
    \begin{subfigure}{0.3\textwidth}
        \includegraphics[width=\textwidth]{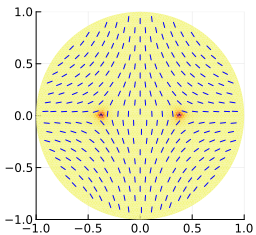}
        \caption{$t = 400$}
    \end{subfigure}
    \begin{subfigure}{0.55\textwidth}
        \includegraphics[width=\textwidth]{img/colorbar.png}
    \end{subfigure}
    \caption{Simulated evolution of a degree -1 tactoid, in which $L_1 = 0.001$ and the initial and boundary conditions are defined by $\mathbf{n}_0 = (-\cos \theta, \sin \theta)$ in (\ref{tactoid_boundary}). Positive eigenvalues and associated eigenvectors, which correspond to the scalar order parameter and the director, are represented by the heat map and quiver plot.}
    \label{fig:tactMinusOne_smallL1}
\end{figure}

\begin{figure}
	\centering
	\begin{subfigure}{0.3\textwidth}
		\includegraphics[width=\textwidth]{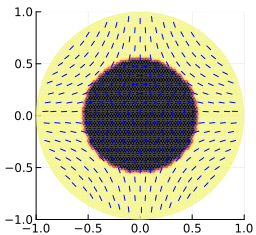}
		\caption{$t = 0.01$}
	\end{subfigure}
	\begin{subfigure}{0.3\textwidth}
		\includegraphics[width=\textwidth]{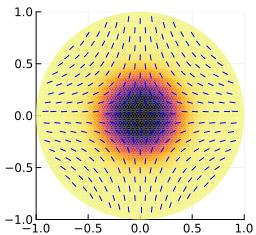}
		\caption{$t = 5$}
	\end{subfigure}
	\begin{subfigure}{0.3\textwidth}
		\includegraphics[width=\textwidth]{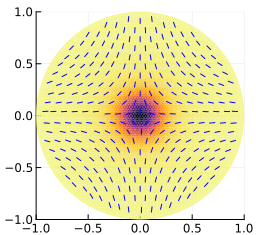}
		\caption{$t = 10$}
	\end{subfigure}
	\begin{subfigure}{0.3\textwidth}
		\includegraphics[width=\textwidth]{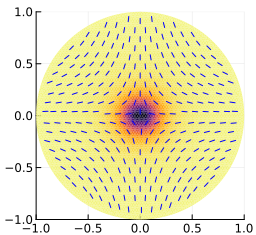}
		\caption{$t = 15$}
	\end{subfigure}
	\begin{subfigure}{0.3\textwidth}
		\includegraphics[width=\textwidth]{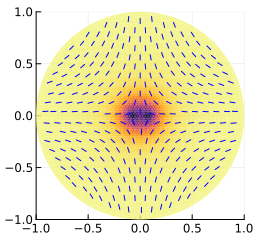}
		\caption{$t = 20$}
	\end{subfigure}
	\begin{subfigure}{0.3\textwidth}
		\includegraphics[width=\textwidth]{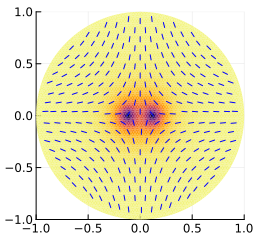}
		\caption{$t = 25$}
	\end{subfigure}
	\begin{subfigure}{0.3\textwidth}
		\includegraphics[width=\textwidth]{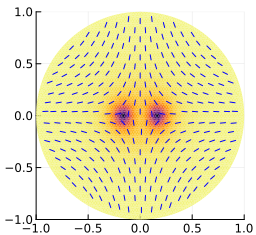}
		\caption{$t = 30$}
	\end{subfigure}
	\begin{subfigure}{0.3\textwidth}
		\includegraphics[width=\textwidth]{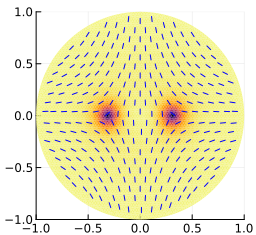}
		\caption{$t = 60$}
	\end{subfigure}
	\begin{subfigure}{0.3\textwidth}
		\includegraphics[width=\textwidth]{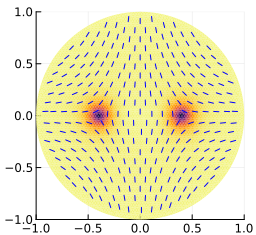}
		\caption{$t = 90$}
	\end{subfigure}
	\begin{subfigure}{0.55\textwidth}
		\includegraphics[width=\textwidth]{img/colorbar.png}
	\end{subfigure}
	\caption{Simulation of~\eqref{eq:oneconstLdG} with $L=0.001$ and the initial and boundary conditions are defined by $\mathbf{n}_0 = (-\cos \theta, \sin \theta)$ in (\ref{tactoid_boundary}). Positive eigenvalues and associated eigenvectors, which correspond to the scalar order parameter and the director, are represented by the heat map and quiver plot.}
	\label{fig:tactMinusOne_LdG}
\end{figure}

In the second trial, we take the director to wind in the reverse direction as the boundary of $\Omega$, by setting
\begin{myeq}
    \mathbf{n}_0 = \begin{pmatrix} -\cos \theta \\ \sin \theta \end{pmatrix}.
\end{myeq}
This leads to a tactoid of topological degree -1. In Figure \ref{fig:tactMinusOne_largeL1}, we observe that when $L_1$ is large, the tactoid shrinks in a radially asymmetrical way until it collapses into a vortex of degree -1, which eventually splits into two defects of degree -1/2 that slowly move apart. Unlike the defects in the previous experiment, these vortices have long V-shaped ``tails'' in which the scalar order parameter is small, and these tails persist as the vortices move. It resembles the experiments in~\cite{golovaty2019phase}. Figure \ref{fig:tactMinusOne_smallL1} shows the tactoid decay when $L_1$ is small, in which event the decay is more radially symmetric but also leaves behind multiple vortices, all but two of which gradually annihilate each other. Finally, in Figure~\ref{fig:tactMinusOne_LdG}, we display the simulation of the standard Landau-de Gennes model~\eqref{eq:oneconstLdG} for $L=0.001$. In contrast to the previous two simulations with the quartic model, in this case the isotropic region shrinks in a spherically symmetric fashion and finally splits into two degree -1/2 defects that move away from each other, just as in the previous two cases. We note that tactoids that have boundaries with kinks, are observed in physical liquid crystal experiments~\cite[Figure 1 and 2]{golovaty2019phase}, thus the behavior the fourth order model displays is desired.

\subsubsection{Degree 0 tactoid}\label{sec:tactZero}
\begin{figure}
    \centering
    \begin{subfigure}{0.3\textwidth}
        \includegraphics[width=\textwidth]{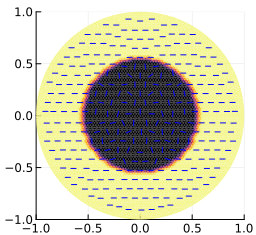}
        \caption{$t = 0.1$}
    \end{subfigure}
    \begin{subfigure}{0.3\textwidth}
        \includegraphics[width=\textwidth]{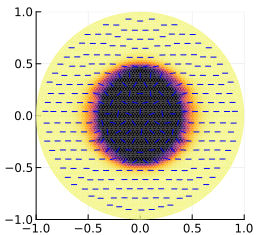}
        \caption{$t = 5$}
    \end{subfigure}
    \begin{subfigure}{0.3\textwidth}
        \includegraphics[width=\textwidth]{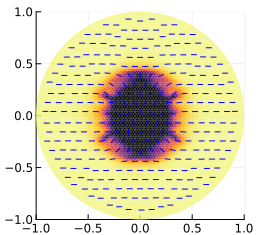}
        \caption{$t = 10$}
    \end{subfigure}
    \begin{subfigure}{0.3\textwidth}
        \includegraphics[width=\textwidth]{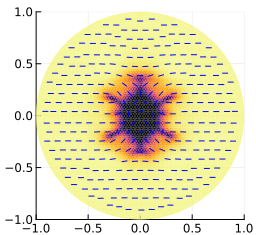}
        \caption{$t = 15$}
    \end{subfigure}
    \begin{subfigure}{0.3\textwidth}
        \includegraphics[width=\textwidth]{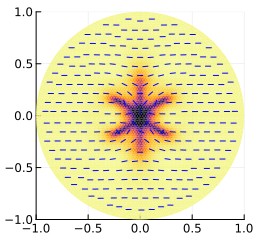}
        \caption{$t = 20$}
    \end{subfigure}
    \begin{subfigure}{0.3\textwidth}
        \includegraphics[width=\textwidth]{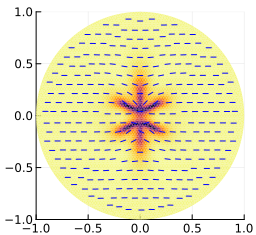}
        \caption{$t = 25$}
    \end{subfigure}
    \begin{subfigure}{0.3\textwidth}
        \includegraphics[width=\textwidth]{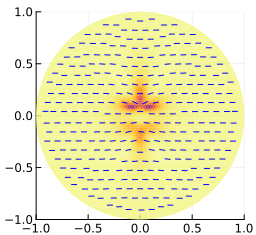}
        \caption{$t = 30$}
    \end{subfigure}
    \begin{subfigure}{0.3\textwidth}
        \includegraphics[width=\textwidth]{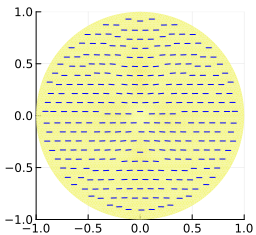}
        \caption{$t = 35$}
    \end{subfigure}
    \begin{subfigure}{0.3\textwidth}
        \includegraphics[width=\textwidth]{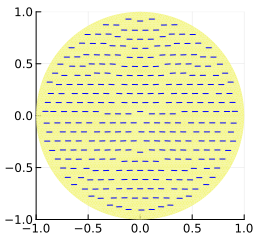}
        \caption{$t = 40$}
    \end{subfigure}
    \begin{subfigure}{0.55\textwidth}
        \includegraphics[width=\textwidth]{img/colorbar.png}
    \end{subfigure}
    \caption{Simulated evolution of a degree 0 tactoid, in which $L_1 = 0.1$ and the initial and boundary conditions are defined by $\mathbf{n}_0 = (1, 0)$ in (\ref{tactoid_boundary}). Positive eigenvalues and associated eigenvectors, which correspond to the scalar order parameter and the director, are represented by the heat map and quiver plot.}
    \label{fig:tactZero_largeL1}
\end{figure}

\begin{figure}
    \centering
    \begin{subfigure}{0.3\textwidth}
        \includegraphics[width=\textwidth]{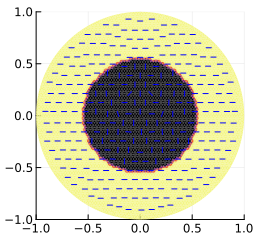}
        \caption{$t = 0.1$}
    \end{subfigure}
    \begin{subfigure}{0.3\textwidth}
        \includegraphics[width=\textwidth]{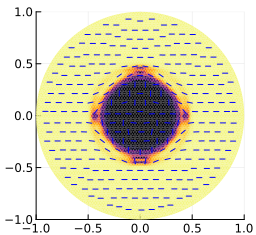}
        \caption{$t = 25$}
    \end{subfigure}
    \begin{subfigure}{0.3\textwidth}
        \includegraphics[width=\textwidth]{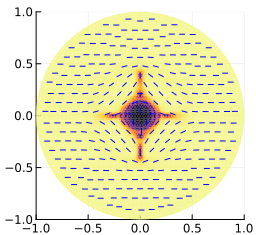}
        \caption{$t = 50$}
    \end{subfigure}
    \begin{subfigure}{0.3\textwidth}
        \includegraphics[width=\textwidth]{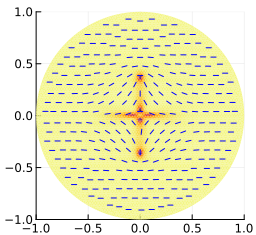}
        \caption{$t = 75$}
    \end{subfigure}
    \begin{subfigure}{0.3\textwidth}
        \includegraphics[width=\textwidth]{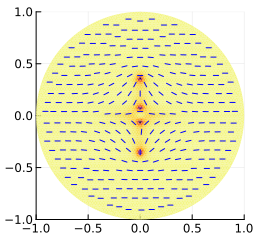}
        \caption{$t = 100$}
    \end{subfigure}
    \begin{subfigure}{0.3\textwidth}
        \includegraphics[width=\textwidth]{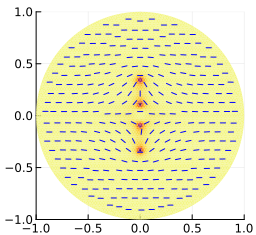}
        \caption{$t = 150$}
    \end{subfigure}
    \begin{subfigure}{0.3\textwidth}
        \includegraphics[width=\textwidth]{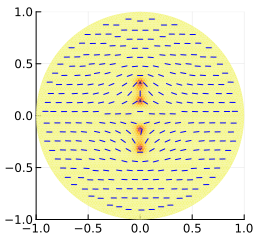}
        \caption{$t = 200$}
    \end{subfigure}
    \begin{subfigure}{0.3\textwidth}
        \includegraphics[width=\textwidth]{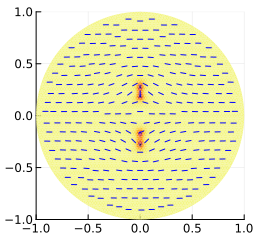}
        \caption{$t = 250$}
    \end{subfigure}
    \begin{subfigure}{0.3\textwidth}
        \includegraphics[width=\textwidth]{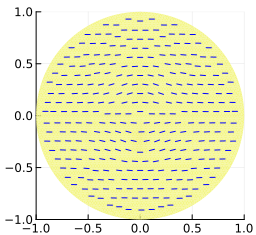}
        \caption{$t = 300$}
    \end{subfigure}
    \begin{subfigure}{0.55\textwidth}
        \includegraphics[width=\textwidth]{img/colorbar.png}
    \end{subfigure}
    \caption{Simulated evolution of a degree 0 tactoid, in which $L_1 = 0.001$ and the initial and boundary conditions are defined by $\mathbf{n}_0 = (1, 0)$ in (\ref{tactoid_boundary}). Positive eigenvalues and associated eigenvectors, which correspond to the scalar order parameter and the director, are represented by the heat map and quiver plot.}
    \label{fig:tactZero_smallL1}
\end{figure}

\begin{figure}
	\centering
	\begin{subfigure}{0.3\textwidth}
		\includegraphics[width=\textwidth]{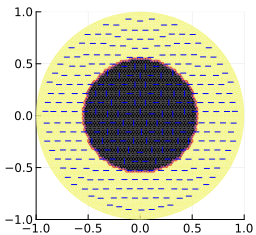}
		\caption{$t = 0.01$}
	\end{subfigure}
	\begin{subfigure}{0.3\textwidth}
		\includegraphics[width=\textwidth]{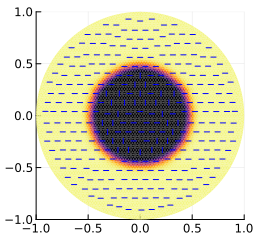}
		\caption{$t = 5$}
	\end{subfigure}
	\begin{subfigure}{0.3\textwidth}
		\includegraphics[width=\textwidth]{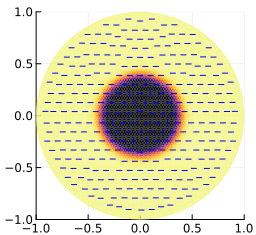}
		\caption{$t = 10$}
	\end{subfigure}
	\begin{subfigure}{0.3\textwidth}
		\includegraphics[width=\textwidth]{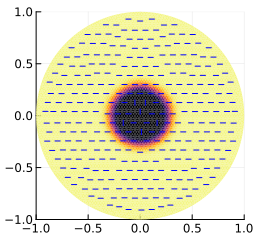}
		\caption{$t = 15$}
	\end{subfigure}
	\begin{subfigure}{0.3\textwidth}
		\includegraphics[width=\textwidth]{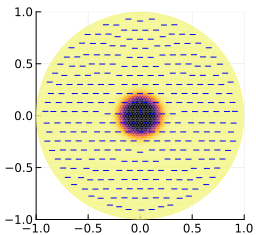}
		\caption{$t = 20$}
	\end{subfigure}
	\begin{subfigure}{0.3\textwidth}
		\includegraphics[width=\textwidth]{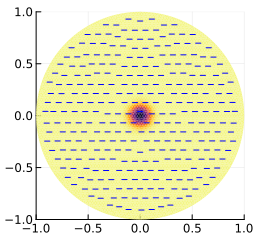}
		\caption{$t = 25$}
	\end{subfigure}
	\begin{subfigure}{0.3\textwidth}
		\includegraphics[width=\textwidth]{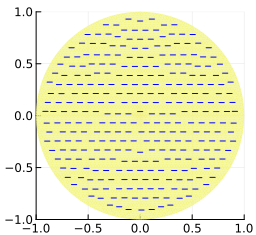}
		\caption{$t = 30$}
	\end{subfigure}
	\begin{subfigure}{0.3\textwidth}
		\includegraphics[width=\textwidth]{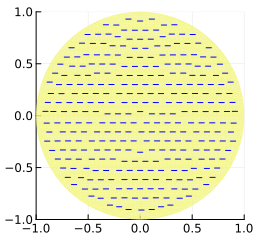}
		\caption{$t = 35$}
	\end{subfigure}
	\begin{subfigure}{0.3\textwidth}
		\includegraphics[width=\textwidth]{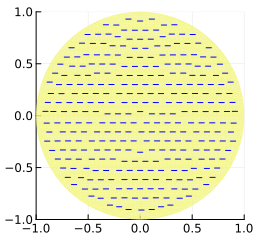}
		\caption{$t = 40$}
	\end{subfigure}
	\begin{subfigure}{0.55\textwidth}
		\includegraphics[width=\textwidth]{img/colorbar.png}
	\end{subfigure}
	\caption{Simulated evolution of a degree 0 tactoid for the quadratic Landau-de Gennes model~\eqref{eq:oneconstLdG} with $L=0.0001$, and the initial and boundary conditions are defined by $\mathbf{n}_0 = (1, 0)$ in (\ref{tactoid_boundary}). Positive eigenvalues and associated eigenvectors, which correspond to the scalar order parameter and the director, are represented by the heat map and quiver plot.}
	\label{fig:tactZero_LdG}
\end{figure}

We now take the constant director
\begin{myeq}
    \mathbf{n}_0 = \begin{pmatrix}
        1 \\ 0
    \end{pmatrix},
\end{myeq}
leading to a tactoid of degree 0. As seen in Figure \ref{fig:tactZero_largeL1}, when $L_1$ is large, this tactoid shrinks until it disappears entirely, leaving behind a uniform sample without vortices. As the tactoid shrinks, six ``arms'' emerge in the shape of an asterisk and persist longer than the surrounding tactoid. Figure \ref{fig:tactZero_smallL1} shows that like the other tactoids, this one also leaves behind vortices when $L_1$ is small; in this case, two with degree $1/2$ and two with degree $-1/2$, which annihilate each other in pairs. In Figure~\ref{fig:tactZero_LdG}, we see a simulation of the standard Landau-de Gennes model~\eqref{eq:oneconstLdG} for $L=0.0001$. In contrast to the quartic model, the isotropic region here shrinks in a spherically symmetric fashion without forming any defects.

\subsubsection{Nematic bubble}\label{sec:bubble}
\begin{figure}
    \centering
    \begin{subfigure}{0.3\textwidth}
        \includegraphics[width=\textwidth]{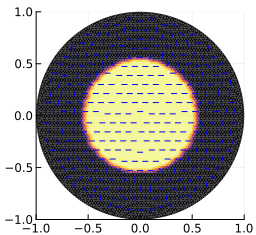}
        \caption{$t = 0.1$}
    \end{subfigure}
    \begin{subfigure}{0.3\textwidth}
        \includegraphics[width=\textwidth]{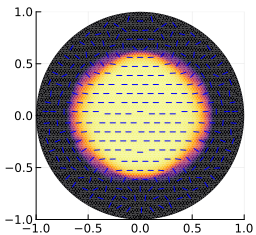}
        \caption{$t = 5$}
    \end{subfigure}
    \begin{subfigure}{0.3\textwidth}
        \includegraphics[width=\textwidth]{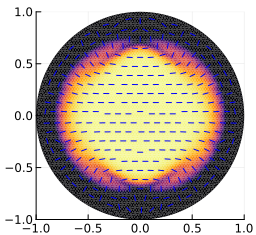}
        \caption{$t = 10$}
    \end{subfigure}
    \begin{subfigure}{0.3\textwidth}
        \includegraphics[width=\textwidth]{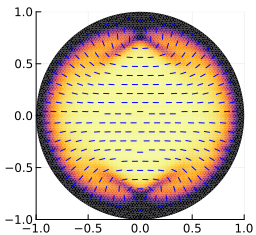}
        \caption{$t = 15$}
    \end{subfigure}
    \begin{subfigure}{0.3\textwidth}
        \includegraphics[width=\textwidth]{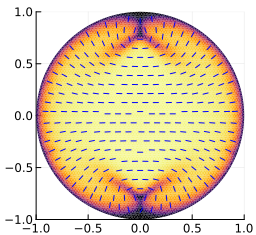}
        \caption{$t = 20$}
    \end{subfigure}
    \begin{subfigure}{0.3\textwidth}
        \includegraphics[width=\textwidth]{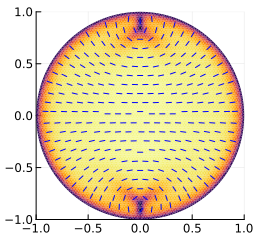}
        \caption{$t = 40$}
    \end{subfigure}
    \begin{subfigure}{0.3\textwidth}
        \includegraphics[width=\textwidth]{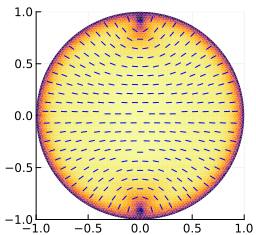}
        \caption{$t = 60$}
    \end{subfigure}
    \begin{subfigure}{0.3\textwidth}
        \includegraphics[width=\textwidth]{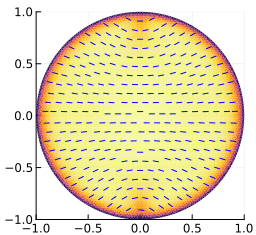}
        \caption{$t = 80$}
    \end{subfigure}
    \begin{subfigure}{0.3\textwidth}
        \includegraphics[width=\textwidth]{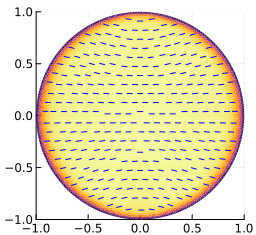}
        \caption{$t = 100$}
    \end{subfigure}
    \begin{subfigure}{0.55\textwidth}
        \includegraphics[width=\textwidth]{img/colorbar.png}
    \end{subfigure}
    \caption{Simulated evolution of a nematic bubble, in which $L_1 = 0.1$ and the initial and boundary conditions are defined by $\mathbf{n}_0 = (1, 0)$ in (\ref{tactoid_boundary}). Positive eigenvalues and associated eigenvectors, which correspond to the scalar order parameter and the director, are represented by the heat map and quiver plot.}
    \label{fig:bubble_largeL1}
\end{figure}

\begin{figure}
    \centering
    \begin{subfigure}{0.3\textwidth}
        \includegraphics[width=\textwidth]{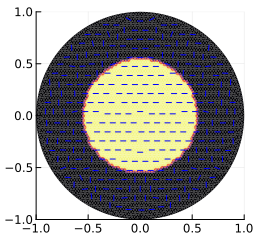}
        \caption{$t = 0.1$}
    \end{subfigure}
    \begin{subfigure}{0.3\textwidth}
        \includegraphics[width=\textwidth]{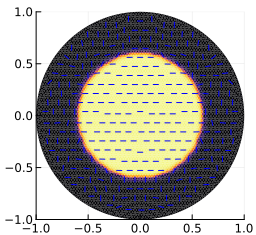}
        \caption{$t = 10$}
    \end{subfigure}
    \begin{subfigure}{0.3\textwidth}
        \includegraphics[width=\textwidth]{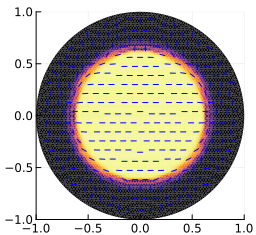}
        \caption{$t = 20$}
    \end{subfigure}
    \begin{subfigure}{0.3\textwidth}
        \includegraphics[width=\textwidth]{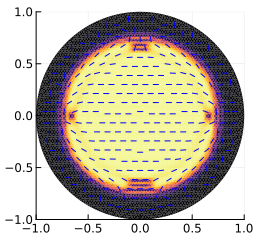}
        \caption{$t = 30$}
    \end{subfigure}
    \begin{subfigure}{0.3\textwidth}
        \includegraphics[width=\textwidth]{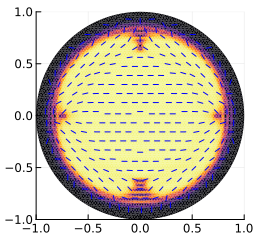}
        \caption{$t = 40$}
    \end{subfigure}
    \begin{subfigure}{0.3\textwidth}
        \includegraphics[width=\textwidth]{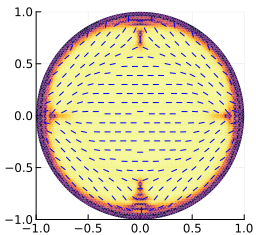}
        \caption{$t = 50$}
    \end{subfigure}
    \begin{subfigure}{0.3\textwidth}
        \includegraphics[width=\textwidth]{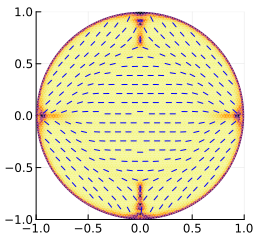}
        \caption{$t = 60$}
    \end{subfigure}
    \begin{subfigure}{0.3\textwidth}
        \includegraphics[width=\textwidth]{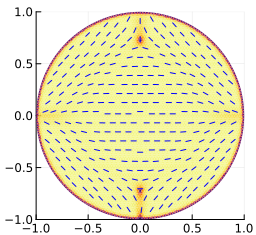}
        \caption{$t = 80$}
    \end{subfigure}
    \begin{subfigure}{0.3\textwidth}
        \includegraphics[width=\textwidth]{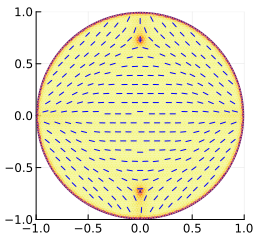}
        \caption{$t = 100$}
    \end{subfigure}
    \begin{subfigure}{0.55\textwidth}
        \includegraphics[width=\textwidth]{img/colorbar.png}
    \end{subfigure}
    \caption{Simulated evolution of a nematic bubble, in which $L_1 = 0.001$ and the initial and boundary conditions are defined by $\mathbf{n}_0 = (1, 0)$ in (\ref{tactoid_boundary}). Positive eigenvalues and associated eigenvectors, which correspond to the scalar order parameter and the director, are represented by the heat map and quiver plot.}
    \label{fig:bubble_smallL1}
\end{figure}

\begin{figure}
	\centering
	\begin{subfigure}{0.3\textwidth}
		\includegraphics[width=\textwidth]{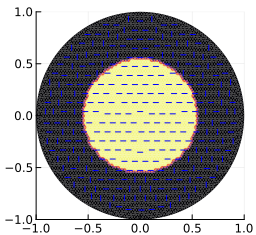}
		\caption{$t = 0.01$}
	\end{subfigure}
	\begin{subfigure}{0.3\textwidth}
		\includegraphics[width=\textwidth]{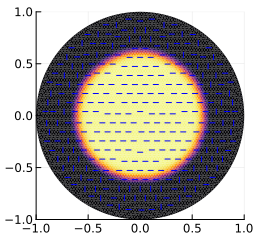}
		\caption{$t = 5$}
	\end{subfigure}
	\begin{subfigure}{0.3\textwidth}
		\includegraphics[width=\textwidth]{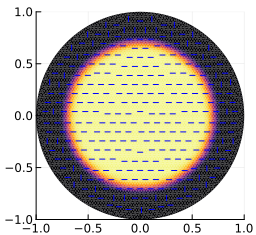}
		\caption{$t = 10$}
	\end{subfigure}
	\begin{subfigure}{0.3\textwidth}
		\includegraphics[width=\textwidth]{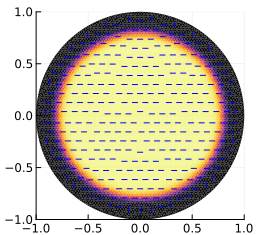}
		\caption{$t = 15$}
	\end{subfigure}
	\begin{subfigure}{0.3\textwidth}
		\includegraphics[width=\textwidth]{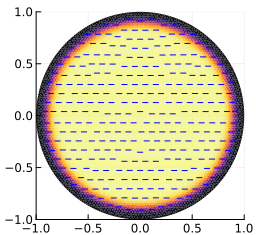}
		\caption{$t = 20$}
	\end{subfigure}
	\begin{subfigure}{0.3\textwidth}
		\includegraphics[width=\textwidth]{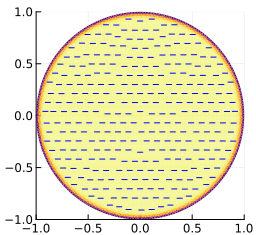}
		\caption{$t = 100$}
	\end{subfigure}
	\begin{subfigure}{0.55\textwidth}
		\includegraphics[width=\textwidth]{img/colorbar.png}
	\end{subfigure}
	\caption{Simulated evolution of a nematic bubble for the standard Landau-de Gennes model~\eqref{eq:oneconstLdG} with $L=0.0001$, in which the initial and boundary conditions are defined by $\mathbf{n}_0 = (1, 0)$ in (\ref{tactoid_boundary}). Positive eigenvalues and associated eigenvectors, which correspond to the scalar order parameter and the director, are represented by the heat map and quiver plot.}
	\label{fig:bubble_LdG}
\end{figure}

\begin{figure}
    \centering
    \includegraphics[width=0.8\linewidth]{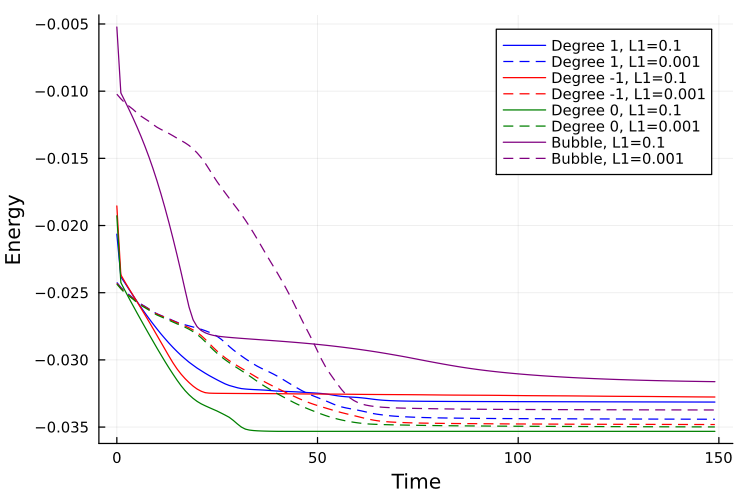}
    \caption{Energies of decaying tactoids and nematic bubbles over time for experiments outlined in Sections \ref{sec:tactOne} (blue), \ref{sec:tactMinusOne} (red), \ref{sec:tactZero} (green), and \ref{sec:bubble} (purple).}
    \label{fig:tactoid_energy}
\end{figure}

Finally, we invert the setup of the previous experiments by including a nematic sample in an isotropic region. The nematic sample initially has constant director
\begin{myeq}
    \mathbf{n}_0 = \begin{pmatrix}
        1 \\ 0
    \end{pmatrix}.
\end{myeq}
We see in both Figure \ref{fig:bubble_largeL1} (in which $L_1$ is large) and Figure \ref{fig:bubble_smallL1} (in which $L_1$ is small), the bubble expands outwards until it reaches the boundary of the sample space. In the former case, two vortices of degree $1/2$ form, but gradually become absorbed into the isotropic ring at the boundary. In the latter case, two defects of degree $-1/2$ form, and these vortices persist within the sample. In Figure~\ref{fig:bubble_LdG}, we show a simulation for the standard quadratic Landau-de Gennes model~\eqref{eq:oneconstLdG}. In this case, no defects form.

\appendix
\section{Tensor calculus lemmas}

\begin{lemma}\label{divgrad3}
    $|\div Q|^2 \leq 3|\nabla Q|^2$ for all $Q \in H^1(\Omega;\RR^{3 \times 3})$.
\end{lemma}
\begin{proof}
    \begin{myeq}
        |\div Q|^2 &= \sum_{i=1}^3 \paren{Q_{i1,1} + Q_{i2,2} + Q_{i3,3}}^2 \\
        &= \sum_{i=1}^3 \paren{3Q_{i1,1}^2 + 3Q_{i2,2}^2 + 3Q_{i3,3}^2 - (Q_{i1,1} - Q_{i2,2})^2 - (Q_{i1,1} - Q_{i3,3})^2 - (Q_{i2,2} - Q_{i3,3})^2} \\
        &\leq \sum_{i=1}^3 \paren{3Q_{i1,1}^2 + 3Q_{i2,2}^2 + 3Q_{i3,3}^2} \leq 3|\nabla Q|^2
    \end{myeq}
\end{proof}

\begin{lemma}\label{curlgrad2}
    $|\curl Q|^2 \leq 2|\nabla Q|^2$ for all $Q \in H^1(\Omega;\RR^{3 \times 3})$.
\end{lemma}
\begin{proof}
    Here we take indices modulo 3 for simplicity; so, for example, $Q_{02,4}$ would mean $Q_{32,1}$.
    \begin{myeq}
        |\curl Q|^2 &= \sum_{i,j=1}^3 (Q_{ji,(i+1)} - Q_{j(i+1),i})^2
        = \sum_{i,j=1}^3 \paren{2Q_{ji,(i+1)}^2 + 2Q_{j(i+1),i}^2 - (Q_{ji,(i+1)} + Q_{j(i+1),i})^2} \\
        &\leq \sum_{i,j=1}^3 \paren{2Q_{ji,(i+1)}^2 + 2Q_{j(i+1),i}^2 + 2Q_{ji,i}^2}
        = 2|\nabla Q|^2
    \end{myeq}
\end{proof}

\begin{lemma}\label{holdermatrix}
    $\myVert{AB}_{L^2} \leq 3 \myVert{A}_{L^\infty} \myVert{B}_{L^2}$ for all $A \in L^\infty(\Omega; \RR^{3 \times 3}), B \in L^2(\Omega; \RR^{3 \times 3})$.
\end{lemma}
\begin{proof}
    \begin{myeq}
        \myVert{AB}_{L^2}^2
        &= \int_\Omega |AB|_F^2 \, dx
        \leq \int_\Omega |A|_F^2 |B|_F^2 \, dx
        = \int_\Omega \paren{\sum_{i,j=1}^3 A_{ij}^2(x)} |B|_F^2 \, dx \\
        &\leq \int_\Omega \paren{9 \max_{i,j,y} A_{ij}^2(y)} |B|_F^2 \, dx
        = 9 \myVert{A}_{L^\infty}^2 \myVert{B}_{L^2}^2
    \end{myeq}
\end{proof}



\bibliographystyle{abbrv}
\bibliography{ref}

\end{document}